\documentclass[11pt,letterpaper]{article}
\usepackage{graphicx}
\usepackage[square,comma,numbers]{natbib}
\usepackage{amsmath}
\usepackage{amssymb}
\usepackage{amsmath}
\usepackage{amsthm}
\usepackage{amsfonts}
\usepackage{mathtools}
\usepackage[super]{nth}
\usepackage{setspace}
\usepackage{enumitem} 
\usepackage{epstopdf}
\usepackage{color}
\usepackage[letterpaper,left=1in,right=1in,top=1in,bottom=1in]{geometry}
\usepackage[ruled,vlined]{algorithm2e}
\usepackage{multirow}
\usepackage{longtable}
\usepackage{rotating}
\usepackage{mathrsfs}
\usepackage{hyperref}

\usepackage{titlesec}
\usepackage{subcaption}
\usepackage{float}
\usepackage{graphicx}
\usepackage{booktabs}
\usepackage{caption}
\usepackage{bm}
\usepackage{tikz}
\usepackage{makecell}
\usepackage[toc,title,page]{appendix}
\usetikzlibrary{shapes,arrows,positioning,calc}

\newtheorem{remark}{Remark}

\titleformat{\paragraph}
{\normalfont\normalsize\bfseries}{\theparagraph}{1em}{}
\titlespacing*{\paragraph}
{0pt}{3.25ex plus 1ex minus .2ex}{1.5ex plus .2ex}

    \title{\Large \bf Graph-Based Imitation and Reinforcement Learning for Efficient Benders Decomposition
    }
\author{
\centerline{\normalsize Bernard T. Agyeman$^{a}$, Zhe Li$^{a}$, Ilias Mitrai$^{b}$\thanks{Corresponding author: Ilias Mitrai. Email: imitrai@che.utexas.edu.}, Prodromos Daoutidis$^{a}$\thanks{Corresponding author: Prodromos Daoutidis. Email: daout001@umn.edu.}}
\vspace{2mm}\\
\centerline{\small $^{a}$Department of Chemical Engineering and Material Science, University of Minnesota,}\\
\centerline{\small Minneapolis, MN 55455, United States.}
\vspace{2mm}\\
\centerline{\small $^{b}$McKetta Department of Chemical Engineering, University of Texas at Austin,}\\
\centerline{\small Austin, TX 78712, United States.}}
\allowdisplaybreaks

\usepackage{url}

\begin{document}
\date{}
\maketitle
\setstretch{1.5}

\begin{abstract}
This work introduces an end-to-end graph-based agent for accelerating the computational efficiency of Benders Decomposition. The agent’s policy is parameterized by a graph neural network which takes as input a bipartite graph representation of the master problem and proposes a candidate solution. The agent is trained using a two-stage approach that combines imitation (IL) and reinforcement learning (RL). IL is used to mimic a solver and obtain a warm-start policy which is then finetuned using RL with a reward signal that balances feasibility and computational efficiency. We augment the agent with a verification mechanism that checks the agent's prediction for feasibility and solution quality. The framework is evaluated in two case studies: (i) an illustrative mixed-integer nonlinear program, where it reduces the solution time by 42\% without loss of solution quality, and (ii) a closed-loop irrigation scheduling problem, where it achieves a 23\% time reduction without compromising water use efficiency.
\end{abstract}

\clearpage
\section{Introduction}

Large-scale mixed-integer optimization problems arise in a multitude of applications in chemical engineering, such as process design, operation, capacity expansion, and parameter estimation \cite{daoutidis2018integrating, pistikopoulos2021process}. A common feature of such problems is the presence of complicating variables, which, if fixed, render the resulting problem easier to solve. This structure is the basis of Benders decomposition, where the original problem is decomposed into a master problem, which is a mixed-integer problem, and a subproblem, which is usually a continuous problem \cite{benders1962partitioning}. The master problem and subproblem are solved sequentially and coordinated via duality, namely, Benders optimality and feasibility cuts. These cuts inform the master problem about the effect of the complicating variables on the subproblem.  Originally developed for mixed-integer linear programming (MILP) problems, Benders decomposition was later extended to convex mixed-integer nonlinear programming (MINLP) problems, giving rise to generalized Benders decomposition (GBD)~\cite{geoffrion1972generalized}.

Despite the wide use of Benders decomposition \cite{conejo2006decomposition,luo2024design,lara2018deterministic,tang2018optimal}, its off-the-shelf implementation is challenging. Specifically, its computational performance depends on the computational complexity of the master problem and subproblem, which are solved iteratively, and the quality of the cuts generated, i.e., the information exchanged between the master problem and subproblem. Over the years, several acceleration techniques have been proposed to improve the computational performance. At a high level, these techniques exploit tailored problem decompositions \cite{magnanti1981accelerating, saharidis2011initialization, crainic2016partial, mitrai2022multicut}, approximate solution of the master problem and subproblem \cite{zakeri2000inexact, larsen2023fast, mitrai2024computationally}, and advanced cut generation and management techniques \cite{magnanti1981accelerating, su2015computational, saharidis2010improving, varelmann2022decoupling}. 

In this paper, we focus on the solution of the master problem, which is an integer optimization problem whose complexity increases as Benders cuts are added during the solution process. Usually, the master problem is solved to global optimality at each iteration, which can be computationally expensive. This has motivated the use of heuristics, such as the variable neighborhood search~\cite{raidl2014speeding}, genetic algorithms~\cite{sohn2012hybrid}, and the tabu search~\cite{jiang2009hybrid}, to quickly obtain feasible solutions. However, these methods generally do not provide optimality guarantees~\cite{raidl2015decomposition}. An alternative is to solve the master problem to local optimality, particularly in the early iterations where the relaxation is weak~\cite{geoffrion1974multicommodity, li2025learning}. Although this locally optimal solution of the master problem can guarantee convergence to the globally optimal solution when the subproblem is convex, it can still lead to high computational effort. 

Given the aforementioned challenges and the iterative nature of Benders decomposition, we propose the application of reinforcement (RL) and imitation (IL) learning for training an agent to solve the master problem at each iteration. In the proposed approach, first expert demonstrations, i.e., master problem solutions using branch and bound, are used to train a policy with imitation learning. In the second step, the policy is fine-tuned with RL to account for the feasibility of the master problem with respect to the predicted binary variables, optimality gap improvement, and solution time of the subproblem. This approach reduces the reliance on mixed-integer branch-and-bound-based optimization solvers and can provide fast approximate solutions in real-time applications.

Machine learning has been used extensively to accelerate the solution of optimization problems ~\cite{bengio2021machine, mitrai2025accelerating} by approximating the optimal solution \cite{bertsimas2022online, ding2020accelerating, gupta2022lookback, paulus2022learning}, tuning monolithic and decomposition-based algorithms \cite{khalil2017learning, morabit2021machine, jia2021benders, mitrai2024takingb, allen2023improvements}, and learning to solve optimization problems \cite{bello2016neural,sonnerat2021learning, bangi2021deep, yoo2021reinforcement}.  Despite these successes, the direct application of RL and IL to solving the master problem poses three challenges. The first is related to the representation of the master problem. Specifically, the structure, i.e., interaction pattern between the variables and constraints in the master problem, changes across iterations since Benders cuts are added. Classical vectorial feature representations of optimization problems cannot capture this evolving interaction pattern \cite{bengio2021machine, mitrai2025accelerating}. The second challenge is related to the lack of feasibility and optimality guarantees when an RL and/or IL agent is used to predict the solution of an optimization problem. Finally, learning to solve combinatorial problems can be challenging due to the large combinatorial space. 

To overcome the aforementioned limitations, we first represent the master problem as a bipartite graph. Graphs have been used extensively to capture the structural interaction between variables and constraints in optimization problems \cite{allman2019decode, jalving2019graph, tang2018optimal}. Moreover, this representation can be augmented with nodal and edge features which capture information about the domain of the variables, right-hand side value of constraints, and coefficient of variables in a constraint. In the context of Benders decomposition, these edge features capture information about the optimality and feasibility of the problem. Given this graph representation, we use Graph Neural Networks (GNNs) to parameterize the agent's policy. GNNs are a class of deep learning architectures designed to operate on graph-structured inputs \cite{bronstein2017geometric}. This approach enables the agent to learn a policy and ultimately predict the values of the complicating variables while accounting for the detailed interaction pattern between the variables and constraints of the master problem at each iteration. This is achieved in two steps, where IL and RL are used to learn and fine-tune the agent's policy. Secondly, we develop a verification mechanism that determines at each iteration whether the agent's predictions should be accepted or if a mixed-integer solver should be used to solve the master problem to global optimality. This verification step relies on criteria such as the prediction confidence of the agent's action, the feasibility of the solution, and the quality of the solution and is the key to guaranteeing convergence of the proposed approach. 

\begin{figure}
    \centering
    \includegraphics[scale=1.0]{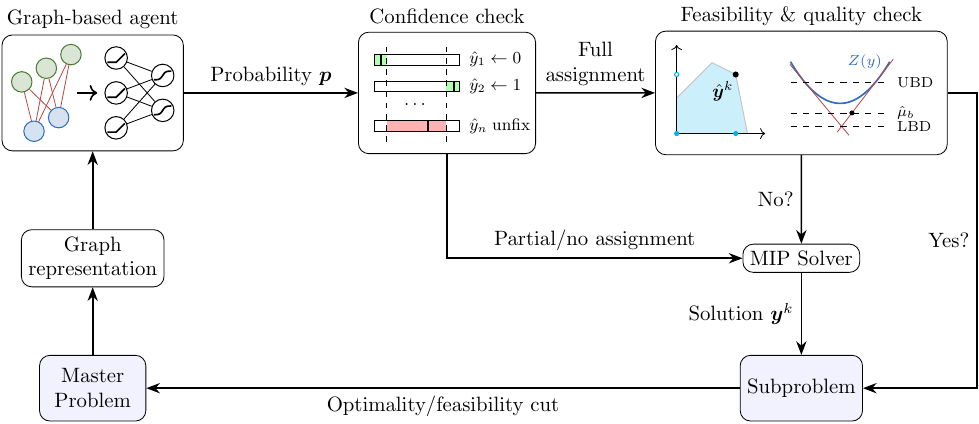}
    \caption{Diagram of the proposed framework.}
    \label{fig:proposed idea}
\end{figure}

Overall, this paper proposes a graph-based hybrid IL / RL agent and verifier for reducing the computational effort related to the Benders decomposition (see Figure~\ref{fig:proposed idea}). Specifically, the contributions of this paper are:
\begin{enumerate}

    \item \textbf{Graph-based agent for master problem acceleration}: A graph-based agent is developed to approximately solve the master problem. The agent’s policy is parameterized by a GNN with a multi-headed output layer, which takes as input a bipartite graph representation of the master problem with nodal and edge features, and outputs a candidate solution. The agent's design combines IL and RL in a complementary way: IL is first used to obtain a warm-start policy by imitating a solver, and this policy initializes the actor of an actor–critic RL framework. RL then fine-tunes the policy using a reward that balances feasibility and computational efficiency.
	
	\item \textbf{Verification mechanism}: A mechanism, called a confidence-based assignment, is developed to coordinate the hybrid use of the agent and the numerical solver within Benders decomposition when solving the master problem. This mechanism ensures that only feasible and sufficiently good solutions are forwarded to the subproblem. It combines prediction confidence, feasibility, and cost consistency checks to determine when to use the agent’s actions and when to defer to the solver.

\end{enumerate}
The proposed framework is evaluated on two optimization problems: (i) an illustrative MINLP, and (ii) a closed-loop irrigation scheduling problem modeled as a mixed-integer model predictive control (MPC).

The remainder of the paper is organized as follows. Section~\ref{sec:preliminaries} introduces the preliminaries, including the class of problems considered, the classical GBD algorithm, graph-based representations optimization problems, and an overview of GNNs, RL and IL. Section~\ref{sec:proposed approach_graph_based} presents the graph-based agent, covering the bipartite graph of the master problem, the imitation learning stage, and the actor–critic RL stage. Section~\ref{sec:proposed_approach_integration} presents the adapted confidence-based approach. Section~\ref{sec:simulations} reports and discusses the results of two case studies, and Section~\ref{sec:conclusion} concludes the paper.

\section{Preliminaries}\label{sec:preliminaries}
\subsection{Problem Formulation}
Without loss of generality, this study focuses on GBD.
We specifically consider the solution of MINLP problems in the following general form:
\begin{equation}
	\begin{aligned}\label{eq:minlp_form}
		\min_{\bm{x},\, \bm{y}} \; &  F(\bm{x}, \bm{y}) \\
		\text{s.t.} \;& \bm{H}(\bm{x}, \bm{y}) = \bm{0}, \\
		& \bm{G}(\bm{x}, \bm{y}) \leq \bm{0}, \\
		& \bm{K}\bm{y} - \bm{b} \leq 0, \\
		& \bm{x} \in \mathcal{X}, \quad \bm{y} \in \{0,1\}^m,
	\end{aligned}
\end{equation}
where \( \bm{x} \in \mathbb{R}^n \) are continuous variables, \( \bm{y} \in \{0,1\}^m \) are binary variables, \( \bm{K} \in \mathbb{R}^{o \times m} \), and \( \bm{b} \in \mathbb{R}^o \). Vector-valued equality and inequality constraints are given by \(\bm{H}:\mathbb{R}^n \times \{0,1\}^m \to \mathbb{R}^p\) and \(\bm{G}:\mathbb{R}^n \times \{0,1\}^m \to \mathbb{R}^q\), respectively.  
The set $\mathcal{X}$ is defined as
\(
\mathcal{X} = \left\{ \bm{x} \in \mathbb{R}^n \,\middle|\, \bm{E}\bm{x} \leq \bm{d},\ \underline{\bm{x}} \leq \bm{x} \leq \overline{\bm{x}} \right\},
\)
where \( \bm{E} \in \mathbb{R}^{l \times n}\), \(\bm{d} \in \mathbb{R}^l\), and \( \underline{\bm{x}},~ \overline{\bm{x}} \) denote the lower and upper bounds on \( \bm{x} \), respectively.  
We focus specifically on a subclass of formulation~\eqref{eq:minlp_form} in which \( \bm{y} \) enters the problem linearly in both the objective function and the constraints, while \( \bm{x} \) appears nonlinearly. Specifically, we assume:
\begin{equation*}
	F(\bm{x}, \bm{y}) \coloneqq f(\bm{x}) + \bm{e}^\top \bm{y}, \quad 
	\bm{H}(\bm{x}, \bm{y}) \coloneqq \bm{h}(\bm{x}) + \bm{A}\bm{y}, \quad 
	\bm{G}(\bm{x}, \bm{y}) \coloneqq \bm{g}(\bm{x}) + \bm{B}\bm{y},
\end{equation*}
where \( \bm{e} \in \mathbb{R}^m \), \( \bm{A} \in \mathbb{R}^{p \times m} \), and \( \bm{B} \in \mathbb{R}^{q \times m} \) are given. The scalar function \( f:\mathbb{R}^n \to \mathbb{R} \), and the vector-valued functions \( \bm{h}:\mathbb{R}^n \to \mathbb{R}^p \) and \( \bm{g}:\mathbb{R}^n \to \mathbb{R}^q \), are assumed to be continuous, differentiable, and convex.

\subsection{Generalized Benders Decomposition (GBD)}\label{sec:classical_gbd}
GBD~\cite{geoffrion1972generalized} proceeds by iteratively refining upper and lower bounds on the optimal solution of the original MINLP problem. For formulation~\eqref{eq:minlp_form}, the upper bound (UBD) is obtained from the subproblem, while the lower bound (LBD) is provided by the master problem.
In iteration \( k \), the subproblem corresponds to the original MINLP with fixed binary variables \( \bm{y}^k \). This problem, denoted as $\mathcal{S}(\bm{y}^k)$, is defined as:
\small
\begin{equation*}
	\begin{aligned}
		Z(\bm{y}^k) = \min_{\bm{x}} \; & F(\bm{x}, \bm{y}^k) \\
		\text{s.t.} \quad & \bm{H}(\bm{x}, \bm{y}^k) = \bm{0}, \\
		& \bm{G}(\bm{x}, \bm{y}^k) \leq \bm{0}, \\
		& \bm{x} \in \mathcal{X}.
	\end{aligned}
	\label{eq:sp_feas}
\end{equation*}
\normalsize
If the subproblem $\mathcal{S}(\bm{y}^k)$ is feasible, its objective value provides an upper bound to the original MINLP. Otherwise, a feasibility subproblem is solved to generate a feasibility cut. This subproblem, denoted as $\mathcal{F}(\bm{y}^k)$, is formulated by introducing slack variables $\bm{\alpha}$ as:
\small
\begin{equation*}
	\begin{aligned}
		\min_{\bm{x},\bm{\alpha}} \; & \sum_{i=1}^{q} \alpha_i \\
		\text{s.t.} \quad & \bm{H}(\bm{x}, \bm{y}^k) = 0, \\
		& G_i(\bm{x}, \bm{y}^k) \leq \alpha_i, \quad i = 1,\dots,q,\\
		&  \alpha_i \geq 0, \quad i = 1,\dots,q,\\		
		& \bm{x} \in \mathcal{X}.
	\end{aligned}
	\label{eq:sp_infeas}
\end{equation*}
Under this decomposition, the original MINLP can be reformulated as 
\begin{equation}
	\begin{aligned}\label{eq:minlp_form_mp}
		\min_{\bm{y}} \; &  Z(\bm{y}) \\
		\text{s.t.} \;& \bm{K}\bm{y} - \bm{b} \leq 0, \\
		& \quad \bm{y} \in \{0,1\}^m.
	\end{aligned}
\end{equation}
Although this is an exact reformulation, the value function of the subproblem $Z(\bm{y})$ is not known explicitly. In Benders decomposition, the value function $Z(\bm{y})$ is approximated via hyperplanes, also known as Benders cuts. For given values of the complicating variables $\bm{y}^{k}$ the Benders optimality cut $\mathcal{C}_o^{(k)}$ is given by:
\begin{equation}\label{eq:optimality_cut}
		\mu_b \geq F(\bm{x}^k, \bm{y}) + \bm{\lambda}_k^\top \bm{H}(\bm{x}^k, \bm{y}) + \bm{\mu}_k^\top \bm{G}(\bm{x}^k, \bm{y}),
\end{equation}
where \( \bm{x}^k \) is the optimal solution to \( \mathcal{S}(\bm{y}^k) \), and \( \bm{\lambda}_k \), \( \bm{\mu}_k \) are the dual multipliers. If the subproblem is infeasible, then the following feasibility cut  $\mathcal{C}_f^{(k)}$ is added:
	\begin{equation}\label{eq:feasibility_cut}
		\overline{\bm{\lambda}}_k^\top \bm{H}(\overline{\bm{x}}^k, \bm{y}) + \overline{\bm{\mu}}_k^\top \bm{G}(\overline{\bm{x}}^k, \bm{y}) \leq 0,
	\end{equation}
where \( \overline{\bm{x}}^k \) is the optimal solution to \( \mathcal{F}(\bm{y}^k) \), and \( \overline{\bm{\lambda}}_k \), \( \overline{\bm{\mu}}_k \) are the dual multipliers.

Overall, the master problem, denoted as $\mathcal{M}$, is equal to
\small
\begin{equation*}
	\begin{aligned}
		\min_{\mu_b,\, \bm{y}} \; & \mu_b \\
		\text{s.t.} \quad 
		& \mu_b \geq F(\bm{x}^k, \bm{y}) + \bm{\lambda}_k^\top \bm{H}(\bm{x}^k, \bm{y}) + \bm{\mu}_k^\top \bm{G}(\bm{x}^k, \bm{y}), \quad \forall k \in K_O \\
		& \overline{\bm{\lambda}}_k^\top \bm{H}(\overline{\bm{x}}^k, \bm{y}) + \overline{\bm{\mu}}_k^\top \bm{G}(\overline{\bm{x}}^k, \bm{y}) \leq 0, \quad \forall k \in K_F \\
		& \bm{K}\bm{y}  - \bm{b} \leq 0,\\
		& \bm{y} \in \{0,1\}^m
	\end{aligned}
	\label{eq:master}
\end{equation*}
\normalsize
where \( K_O \) and \( K_F \) denote the sets of iterations where the subproblem is feasible and infeasible, respectively. Since the number of Benders feasibility and optimality cuts can be very large, they are added adaptively during the solution process. Thus, at each iteration, the solution of the master problem provides a lower bound (LBD). The original MINLP constraints involving only the binary variables, which we refer to as pure binary constraints,  are denoted as $\mathcal{C}_p$:
	\begin{equation} \label{eq:pure_binary_constraints}
		\bm{K}\bm{y}  - \bm{b} \leq 0.
	\end{equation}

\noindent
The complete GBD algorithm, described in Algorithm~\ref{alg:gbd_alg}, proceeds by alternating between solving the master problem and subproblem until the gap between UBD and LBD falls below a specified tolerance \( \varepsilon \). 

\begin{algorithm}[t]
	\caption{Generalized Benders decomposition~\cite{geoffrion1972generalized}.}
	\footnotesize
	\label{alg:gbd_alg}
	\DontPrintSemicolon
	\SetAlgoLined
	\KwIn{\( \bm{y}^0 \in \{0,1\}^m \), \( \varepsilon > 0 \)}
	\textbf{Initialize: } \( k \gets 0 \), LBD \( \gets -\infty \), UBD \( \gets +\infty \), \( K_O \gets \emptyset \), \( K_F \gets \emptyset \)\;
	
	\While{\( \text{UBD} - \text{LBD} > \varepsilon \)}{
		Solve subproblem \( \mathcal{S}(\bm{y}^k) \)\;
		
		\eIf{subproblem \( \mathcal{S}(\bm{y}^k) \) is feasible}{
			Add optimality cut 
			\(
			\mu_b \geq F(\bm{x}^k, \bm{y}) 
			+ \bm{\lambda}_k^\top \bm{H}(\bm{x}^k, \bm{y}) 
			+ \bm{\mu}_k^\top \bm{G}(\bm{x}^k, \bm{y}) \) to master problem \( \mathcal{M} \) \;
			\( \text{UBD} \gets \min(\text{UBD}, Z(\bm{y}^k)) \)\;
			\( K_O \gets K_O \cup \{k\} \)\;
		}{
			Solve subproblem \( \mathcal{F}(\bm{y}^k) \)\;
			Add feasibility cut 
			\( \;
			\overline{\bm{\lambda}}_k^\top \bm{H}(\overline{\bm{x}}^k, \bm{y}) 
			+ \overline{\bm{\mu}}_k^\top \bm{G}(\overline{\bm{x}}^k, \bm{y}) 
			\leq 0 \) to master problem \( \mathcal{M} \) \;
			\( K_F \gets K_F \cup \{k\} \)\;
		}
		
		Solve master problem \( \mathcal{M} \) to obtain \( \bm{y}^{k+1} \), \( \mu_b^k \)\;
		\( \bm{y}^{k} \gets \bm{y}^{k+1} \)\;
		LBD \( \gets \mu_b^k \)\; 
		\( k \gets k + 1 \)\;
	}
\end{algorithm}

\subsection{Graph Representation of Optimization Problems}\label{sec:graph_rep_opt}
Graphs $\mathcal{G}$ are mathematical structures consisting of nodes $\mathcal{V}$ and edges $\mathcal{E}$. Bipartite graphs are especially useful for representing optimization problems~\cite{gasse2019exact, tang2018optimal, allman2019decode}. A bipartite graph is one whose nodes can be partitioned into two disjoint subsets such that each edge connects a node from one subset to a node in the other. 
In the context of optimization problems, one subset of nodes corresponds to the variables, while the other corresponds to the constraints. An edge is introduced between a variable node and a constraint node whenever the variable appears in that constraint \cite{allman2019decode}.  The connectivity between the nodes of the graph is captured by the adjacency matrix ${A}_d \in \{0,1\}^{|\mathcal{V}|\times|\mathcal{V}|}$, where $({A}_d)_{ij} = 1$ if an edge exists between nodes $i$ and $j$, and $0$ otherwise.

The nodes (variable and constraint) and the edges of the bipartite graph can be augmented with features that encode additional information about the optimization problem. These node features are collected in a node feature matrix ${X}_f$. Edges may include non-zero coefficients of the constraint matrix as features, which are stored in an edge feature matrix $X_{e}$. Together, the matrices ${A}_d$, ${X}_f$, and $X_{e}$ provide a compact way to represent the bipartite graph of an optimization problem.

\subsection{Graph Neural Networks}\label{sec:gnn_prelims}
Graph Neural Networks (GNNs) are a class of neural networks designed to operate on graph-structured data. Given a graph \( \mathcal{G} = (\mathcal{V}, \mathcal{E}) \), a GNN learns permutation-invariant functions that capture both local and global structure. Most GNN architectures follow a message-passing paradigm, where each node updates its representation by aggregating information from its neighbors. Common examples include the Graph Convolutional Network (GCN)~\cite{kipf2016semi} and the Edge-Conditioned Convolution (ECC)~\cite{simonovsky2017dynamic}.

Of particular relevance to this work is the ECC, which extends GNNs to incorporate edge features into the message-passing process. Let \( X_f^{(\ell)} \) denote the node feature matrix at layer \( \ell \), \( X_e \) the edge feature matrix, and \( \mathcal{N}(i) \) the set of neighbors of node \( i \) as defined by the adjacency matrix \( A_d \). The ECC update rule is expressed as~\cite{simonovsky2017dynamic}:
\begin{align}
    X_f^{(0)} &= X_f,\\
    X_f^{(\ell)}(i) &= \frac{1}{|\mathcal{N}(i)|} \sum_{j \in \mathcal{N}(i)} 
    \mathcal{F}^{(\ell)}(X_e(j,i); \bm{w}^{(\ell)}) \, X_f^{(\ell-1)}(j) + \bm{b}^{(\ell)},
    \quad \ell = 1, \dots, L,
    \label{eq:ecc}
\end{align}
where \( X_f^{(\ell)} \) denotes the node embedding at layer \( \ell \), \( \mathcal{F}^{(\ell)} \) is a learnable function typically implemented as a multi-layer perceptron (MLP) parameterized by weights \( \bm{w}^{(\ell)} \), and \( \bm{b}^{(\ell)} \) is a bias term.

The node embeddings at the final layer, \( X_f^{(L)} \), can serve as inputs to another MLP that produces task-specific outputs. Depending on the learning objective, GNNs can perform tasks (e.g., classification, regression) at the node level, edge level, or graph level. For graph-level tasks, pooling and readout operations are combined with the GNN layers to aggregate the node embeddings into a fixed-size graph representation. An important property of GNN architectures is their ability to process graphs of different sizes and topologies using the same set of learned parameters. This makes GNNs suitable for applications involving evolving graph structures, such as the bipartite graph representation of the master problem in GBD, where the addition of cuts at each iteration yields a graph whose structure evolves across iterations.

\subsection{Reinforcement Learning (RL)}
RL is a machine learning paradigm in which an agent learns to make decisions by interacting with an environment and seeks to maximize cumulative rewards over time.  The standard formalism for RL is the Markov Decision Process (MDP), which provides a mathematical framework for modeling sequential decision-making problems. A finite MDP is defined by a tuple $(\mathcal{S}, \mathcal{A}, \mathcal{P}, \mathcal{R})$, where $\mathcal{S}$ is a finite set of states, $\mathcal{A}$ is a finite set of actions, $\mathcal{P}(s' \mid s,a)$ defines the probability of moving to state $s'$ and receiving a reward $r$ when action $a$ is taken in state $s$, and $\mathcal{R}$ is the reward function. 
At each time step $t$, the agent observes a state $s_t \in \mathcal{S}$, selects an action $a_t \in \mathcal{A}$, and receives a scalar reward $r_t$, while the environment transitions to a new state $s_{t+1}$ according to $\mathcal{P}$. The agent’s objective is to learn a policy $\pi: \mathcal{S} \rightarrow \mathcal{A}$ that maximizes the expected cumulative discounted reward
\(
\mathbb{E} \left[ \sum_{t=0}^{\infty} \gamma^t r_t \right],
\)
where $\gamma \in [0,1)$ is the discount factor.

RL methods are generally classified into value-based and policy-based methods. Value-based methods aim to learn a value function that estimates the expected return under a given policy. The optimal policy is then derived indirectly by selecting actions that maximize the learned value function. Policy-based methods, in contrast, directly optimize a parameterized policy via gradient ascent on the expected return. Central to this work are actor–critic methods, such as Proximal Policy Optimization (PPO)~\cite{schulman2017proximal}, which combine the strengths of value-based and policy-based paradigms. These methods maintain two distinct models: an actor, which represents a parameterized policy, and a critic, which estimates a value function to guide the actor’s updates.

\subsection{Imitation Learning (IL)}
IL is a learning paradigm in which an agent seeks to learn a policy by mimicking expert behavior. 
Formally, the IL problem shares structural similarities with RL. However, the key distinction lies in the absence of an explicit reward signal. Rather, the agent is provided with a collection of expert trajectories \( \xi \), where each trajectory \( \xi = \{(s_0, a_0), (s_1, a_1), \dots \} \) comprises state-action pairs sampled from an expert policy \( \pi^* \). The objective of imitation learning is to recover a policy \( \pi \) that approximates the behavior of the expert.

Relevant to this work is behavioral cloning, which formulates the imitation learning task as a supervised learning problem.
Behavioral cloning is simple to implement and effective in settings with sufficient and representative expert data. In behavioral cloning, the agent is trained to minimize the discrepancy between its predicted actions and those of the expert. Given a data set of expert demonstrations \(\Xi\), the objective is to solve:
\[
\hat{\pi}^* = \arg\min_{\pi} \sum_{\xi \in \Xi} \sum_{s \in \xi} \mathcal{L}(\pi(s), \pi^*(s)),
\]
where \( \mathcal{L} \) is a suitable loss function, \( \pi^*(s) \) denotes the action of the expert in the state \( s \), and \( \hat{\pi}^* \) is the learned policy that imitates the expert.

\section{Graph-Based Agent}\label{sec:proposed approach_graph_based}
We present the proposed graph-based agent to solve the master problem. The agent takes a bipartite graph representation of the master problem as its input and predicts the candidate values of the binary variables in the master problem. The development of the agent follows a two-stage approach, which integrates IL and RL. In the first stage, the agent is trained through IL to mimic a mixed-integer programming (MIP) solver.  The IL agent policy then serves as the initial policy within an RL framework. In the second stage, this policy is fine-tuned through interactions with the GBD environment. 

The remainder of this section is organized as follows: first, we present the bipartite graph of the master problem; second, we describe the architecture and training procedure employed in the IL stage; and finally, we outline the MDP formulation used for RL training.

\begin{figure}[t]
	\centerline{\includegraphics[width=\textwidth]{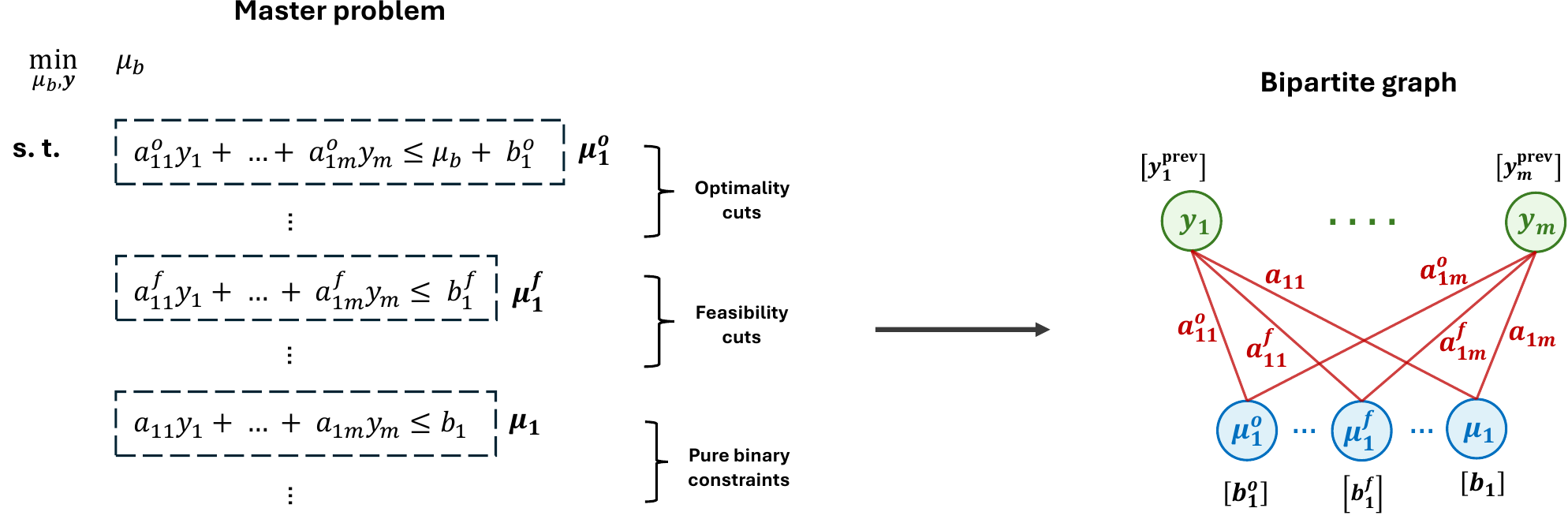}}
	\caption{Bipartite graph of the master problem. The constraints are linear because the pure binary constraints are linear and the binary variables enter  the objective and constraints of the original MINLP linearly.}
	\label{fig:graph_mp_example}
\end{figure}

\subsection{Bipartite Graph Representation of the Master Problem}\label{sec:bipartite_rep_of_mp}
Each instance of the master problem is represented as a bipartite graph, following the general framework introduced in Section~\ref{sec:graph_rep_opt}.
Figure~\ref{fig:graph_mp_example} illustrates a generic master problem instance and its graph encoding; the key components of the graph are summarized below.

\begin{enumerate}[label=(\roman*)]
	\item \textbf{Variable nodes}: Each variable node corresponds to a binary variable in the master problem. The continuous variable $\mu_b$, which appears in the optimality cuts, is not included as a variable node. A justification for this design choice is given in Remark~\ref{rmk:mu_exclusion}.
	
	\item \textbf{Variable node features}: 
    The value of each binary variable from the previous master problem iteration was used as the variable node feature. This choice is consistent with related studies, such as~\cite{gasse2019exact}. 
    In Figure~\ref{fig:graph_mp_example}, these features are denoted as \( \left[y_1^{\text{prev}}, \dots, y_m^{\text{prev}}\right] \).
	
	\item \textbf{Constraint nodes}: Each constraint node represents a constraint in the master problem, including feasibility cuts ($\mu_{(\cdot)}^f$), optimality cuts ($\mu_{(\cdot)}^o$), and pure binary constraints ($\mu_{(\cdot)}$), as illustrated in Figure~\ref{fig:graph_mp_example}.

    \item \textbf{Constraint node features}:
The constant term on the right-hand side (RHS) of each constraint (i.e., the optimality cuts, feasibility cuts, and pure-binary constraints) serves as the feature of the corresponding constraint node. In Figure~\ref{fig:graph_mp_example}, the constraint-node features of the optimality cuts, feasibility cuts, and pure-binary constraints are denoted as $b_{k}^{o}$ for $k \in K_O$, $b_{k}^{f}$ for $k \in K_F$, and $b_{t}$ for $t \in [1,\dots,s]$, respectively. 
Given formulation~\eqref{eq:minlp_form}, as well as the expressions for the optimality cuts (inequality~\ref{eq:optimality_cut}), feasibility cuts (inequality~\eqref{eq:feasibility_cut}), and pure-binary constraints (inequality~\eqref{eq:pure_binary_constraints}), the various constraint-node features can be obtained as follows. 

For the \emph{optimality cuts},
\[
\mu_b \ge f(\bm{x}^k) + \bm{e}^\top\bm{y}
+ \bm{\lambda}_k^\top\!\big(\bm{h}(\bm{x}^k)+\bm{A}\bm{y}\big)
+ \bm{\mu}_k^\top\!\big(\bm{g}(\bm{x}^k)+\bm{B}\bm{y}\big),
\qquad \forall\, k \in K_O,
\]
the corresponding constraint-node feature is
\[
b_k^{o} = -\!\left[f(\bm{x}^k) + \bm{\lambda}_k^\top\bm{h}(\bm{x}^k) + \bm{\mu}_k^\top\bm{g}(\bm{x}^k)\right],
\qquad \forall\, k \in K_O,
\]
where $\bm{x}^k$ is the optimal solution of the subproblem $\mathcal{S}(\bm{y}^k)$ when it is feasible, and $\bm{\lambda}_k$ and $\bm{\mu}_k$ are the dual multipliers.

For the \emph{feasibility cuts},
\[
\overline{\bm{\lambda}}_k^\top\!\big(\bm{h}(\overline{\bm{x}}^k)+\bm{A}\bm{y}\big)
+ \overline{\bm{\mu}}_k^\top\!\big(\bm{g}(\overline{\bm{x}}^k)+\bm{B}\bm{y}\big) \le 0,
\qquad \forall\, k \in K_F,
\]
the corresponding constraint-node feature is
\[
b_k^{f} = -\!\left[\overline{\bm{\lambda}}_k^\top\bm{h}(\overline{\bm{x}}^k) + \overline{\bm{\mu}}_k^\top\bm{g}(\overline{\bm{x}}^k)\right],
\qquad \forall\, k \in K_F,
\]
where $\overline{\bm{x}}^k$ is the optimal solution of the feasibility subproblem $\mathcal{F}(\bm{y}^k)$, and $\overline{\bm{\lambda}}_k$ and $\overline{\bm{\mu}}_k$ are the dual multipliers.

For the \emph{pure-binary constraints},
the constraint-node features correspond directly to the entries of the vector~$\bm{b}$ in inequality~\eqref{eq:pure_binary_constraints}. 

Collectively, the constraint-node features are 
\(\{b_k^{o}\}_{k\in K_O}\), 
\(\{b_k^{f}\}_{k\in K_F}\), and 
\(\{b_t\}_{t=1}^{s}\),
which together form the constraint-feature rows in \(X_f\).
	
	\item \textbf{Edges}: An edge exists between a variable node and a constraint node if the corresponding binary variable appears in that constraint.  Connectivity among nodes is captured by the adjacency matrix $A_d$.

\item \textbf{Edge features}:
Each edge feature corresponds to the nonzero coefficient of a binary variable in the associated constraint. In Figure~\ref{fig:graph_mp_example}, the edge features of the optimality cuts, feasibility cuts, and pure-binary constraints are denoted as $a_{kj}^{o}$, $a_{kj}^{f}$, and $a_{tj}$, respectively. 
Given formulation~\eqref{eq:minlp_form}, as well as the expressions of the optimality cuts, feasibility cuts, and pure-binary constraints, the various edge features are derived as follows.

For the \emph{optimality cuts},
\[
\mu_b \ge f(\bm{x}^k) + \bm{e}^\top\bm{y}
+ \bm{\lambda}_k^\top\!\big(\bm{h}(\bm{x}^k)+\bm{A}\bm{y}\big)
+ \bm{\mu}_k^\top\!\big(\bm{g}(\bm{x}^k)+\bm{B}\bm{y}\big),
\qquad \forall\, k \in K_O,
\]
the edge feature linking variable \(y_j\) to the \(k\)th optimality cut is
\[
a_{kj}^{o} = e_j + (\bm{A}^\top\bm{\lambda}_k)_j + (\bm{B}^\top\bm{\mu}_k)_j,
\qquad \forall\, j \in [1,\dots,m],\; k \in K_O.
\]

For the \emph{feasibility cuts},
\[
\overline{\bm{\lambda}}_k^\top\!\big(\bm{h}(\overline{\bm{x}}^k)+\bm{A}\bm{y}\big)
+ \overline{\bm{\mu}}_k^\top\!\big(\bm{g}(\overline{\bm{x}}^k)+\bm{B}\bm{y}\big) \le 0,
\qquad \forall\, k \in K_F,
\]
the edge feature linking variable \(y_j\) to the \(k\)th feasibility cut is
\[
a_{kj}^{f} = (\bm{A}^\top\overline{\bm{\lambda}}_k)_j + (\bm{B}^\top\overline{\bm{\mu}}_k)_j,
\qquad \forall\, j \in [1,\dots,m],\; k \in K_F.
\]

For the \emph{pure-binary constraints},
the edge features correspond to the entries of $\bm{K}$ in inequality~\ref{eq:pure_binary_constraints}, that is,
\[
a_{tj} = K_{tj},
\qquad \forall\, j \in [1,\dots,m],\; t \in [1,\dots,s].
\]

Collectively, the edge features are 
\(\{a_{kj}^{o}\}_{k\in K_O}\), 
\(\{a_{kj}^{f}\}_{k\in K_F}\), and 
\(\{a_{tj}\}_{t=1}^{s}\) for all \(j \in [1,\dots,m]\),
which together form the edge-feature matrix \(X_e\).

\end{enumerate}
Under this graph encoding, the number of variable nodes remains fixed across master problem instances, while the number of constraint nodes and edges grows over iterations as cuts are added. Consequently, the matrices $A_d$, $X_f$, and \(X_e\) grow in size as the GBD algorithm progresses. Furthermore, to enhance training stability and predictive performance, the entries of matrices $X_f$ and \(X_e\) are normalized.
\begin{remark}\label{rmk:mu_exclusion}
	The continuous variable \( \mu_b \) appears in the optimality cuts in the master problem, which could justify its inclusion as a variable node in the graph. We evaluated this setting in our experiments but found that including \( \mu_b \) in the graph consistently reduced the predictive accuracy of the imitation learning agent. Based on this empirical evidence, \( \mu_b \) is excluded from the master problem graph.
\end{remark}

\subsection{IL Stage}\label{sec:il_stage_proposed_approach}
\begin{figure}[t]
	\centerline{\includegraphics[width=\textwidth]{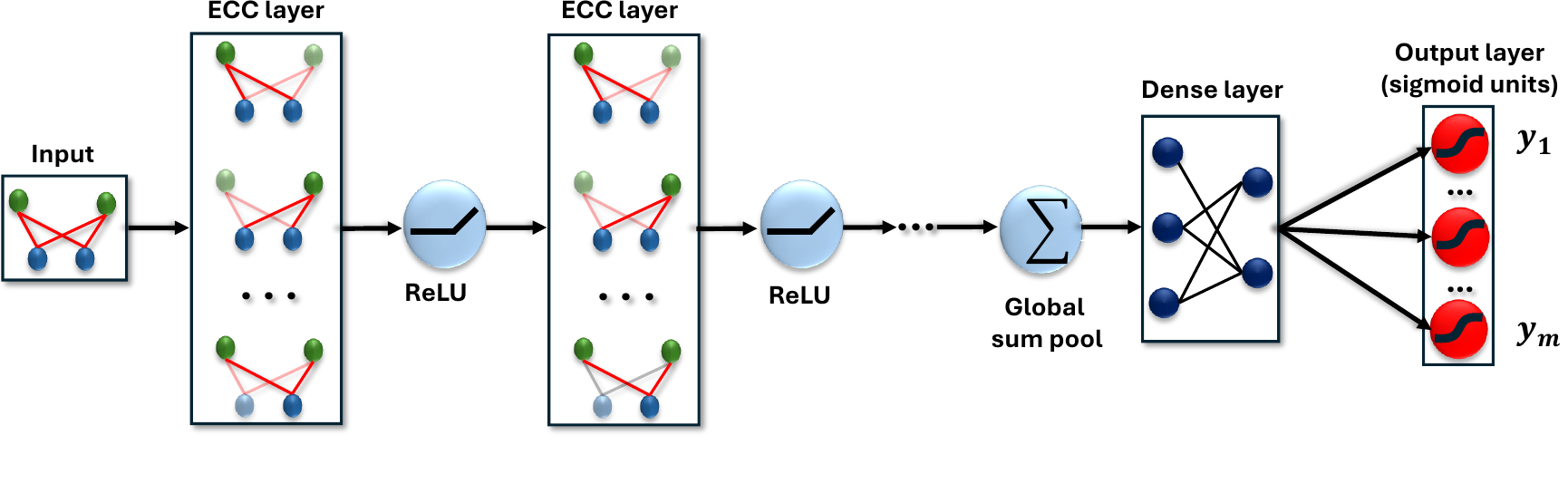}}
	\caption{A schematic of the GNN with multi-headed sigmoid output layer.}
	\label{fig:gnn_diagram}
\end{figure}
The agent is parameterized by a GNN to handle the graph-structured nature of its inputs and to exploit the effectiveness of deep neural networks in approximating policies. The learning problem at this stage is formulated as a graph-level multi-label classification task. A schematic of the agent is shown in Figure~\ref{fig:gnn_diagram}, and its key components are described below:
\begin{itemize}
	\item \textbf{Input layer}: The bipartite graph of the master problem.
	\item \textbf{Graph layer(s)}: One or more ECC layers.
	\item \textbf{Readout layer}:  A global sum pooling operation that aggregates node embeddings to produce a graph-level representation.
	\item \textbf{Dense layer(s)}: Fully connected layers are applied after pooling to increase the network's capacity to capture complex interactions.
	\item \textbf{Output layer}: A layer with as many units as binary variables in the master problem. Each unit uses a sigmoid activation function to produce a value in the range \([0, 1]\). This design yields a multi-headed output architecture, where each binary variable is predicted independently.
\end{itemize}
Following the specification of the agent's architecture, it is trained using behavioral cloning. During training, the binary cross-entropy loss function is employed for the multi-label classification task.
For a given set of ground truth labels \( \bm{y} \)  and predicted values \( \hat{\bm{y}} \), the cross-entropy loss function \( \mathcal{L} (\cdot) \) is defined as:
\[
\mathcal{L}(\bm{y},\hat{\bm{y}}) \coloneqq \frac{1}{m} \sum_{i=1}^{m} \left[ -y_i \log \hat{y}_i - (1 - y_i)\log(1 - \hat{y}_i) \right].
\]

\begin{remark}\label{rmk:scalability}
	An alternative to the multi-headed output structure is to model the task as a discrete action selection problem, where each possible combination of binary assignments corresponds to a unique action. For example, with 3 binary variables, there are \(2^3 = 8\) possible assignments. Although conceptually simple, this approach scales poorly as the number of actions grows exponentially with the number of binary variables. By contrast, the adopted multi-headed architecture scales linearly with the number of decision variables and is therefore more suitable for large-scale problems. In addition, the units of the output layer share a common graph, thereby ensuring that the predictions are informed by the global structure of the master problem.
\end{remark}

\subsection{RL Stage}
In this stage, the agent resulting from the first stage is fine-tuned within an actor–critic RL framework. The actor is initialized with the imitation learning policy and therefore adopts the same architecture as the imitation learning agent. To handle the graph-structured nature of its inputs and to exploit the effectiveness of deep neural networks for function approximation, the critic is also implemented as a GNN. Its architecture follows the same design as that of the actor (see Section~\ref{sec:il_stage_proposed_approach}) but differs in the output layer. Unlike the actor, which possesses a multi-headed output layer with sigmoid activation, the critic has a single-unit output layer without activation to approximate the state-value function. The key components of the MDP formulation that govern the agent’s decision-making process are described below.

\subsubsection{State}  
The state of the agent is the bipartite graph of the master problem instance, which is compactly represented by the adjacency matrix $A_d$, the node feature matrix $X_f$, and the edge feature matrix \(X_e\). Together, $A_d$, $X_f$, and \(X_e\) define the state and serve as input to the actors and critic networks.

\subsubsection{Action}\label{sec:action}  
The action $\bm{a}$ of the agent belongs to the discrete space $\{0,1\}^m$, which corresponds to the $m$ binary decision variables of the master problem. Each action therefore corresponds to a complete assignment of values to these binary decision variables. The design of the actor, particularly its output layer, makes the policy inherently stochastic. Each unit of the actor’s output layer is equipped with a sigmoid activation, producing a value in the interval \([0,1]\). This value is interpreted as the probability of assigning 1 to the corresponding binary variable. Thus, each output $p_i$ defines a Bernoulli distribution, from which the binary decision can be sampled. The agent’s action $\bm{a}$ is obtained as a vector of independent samples from the Bernoulli distributions parameterized by $\bm{p}$. Given the policy $\pi_{\theta}$, the action selection proceeds as follows:
\begin{equation} \label{eq:action_selection}
	\begin{aligned}
		\pi_{\theta}(A_d, X_f, L) &\coloneqq \bm{p} \coloneqq (p_1,\dots,p_m), \quad p_i \in [0,1], \\
		a_i &\sim \text{Bernoulli}(p_i),~i = 1,\dots,m, \\
		\bm{a} &\coloneqq (a_1,\dots,a_m) \in \{0,1\}^m .
	\end{aligned}
\end{equation}

\subsubsection{Transition Dynamics}  
In iteration $k$, the current agent policy uses the matrices $A_{d_k}$, $X_{f_k}$, and $X_{e_k}$ to produce action $\bm{a}_k$ according to Equation~\ref{eq:action_selection}. However, $\bm{a}_k$ cannot be passed directly to the subproblem, as its feasibility with respect to the master problem is not guaranteed. Thus, $\bm{a}_k$ is first checked for feasibility. This check involves evaluating feasibility cuts (that is, the left-hand side (LHS) of inequality~\ref{eq:feasibility_cut}) and the pure binary constraints (that is, the LHS of inequality~\ref{eq:pure_binary_constraints}) using $\bm{a}_k$. If these constraints are satisfied, $\bm{a}_k$ is passed to the subproblem, which is solved to generate the appropriate cut. Conversely, if $\bm{a}_k$ is infeasible, the agent receives a penalty, and an MIP solver is invoked to solve the master problem. The resulting solution, $\bm{\overline{y}}_k$, is then passed to the subproblem. In either case, a transition occurs in which the master problem is updated by adding a cut, which modifies the matrices $A_d$, $X_f$, and $X_e$ to yield $A_{d_{k+1}}$, $X_{f_{k+1}}$, and $X_{e_{k+1}}$. In summary, under $\bm{a}_k$ (or $\bm{\overline{y}}_k$), the state of the master problem transitions from $(A_{d_k}, X_{f_k}, X_{e_k})$ to $(A_{d_{k+1}}, X_{f_{k+1}}, X_{e_{k+1}})$, which defines the dynamics of the agent’s environment.

\subsubsection{Reward Design}
The agent's objective is to choose actions that (i) render the master problem feasible, (ii) improve the GBD bound gap, and (iii) reduce the solution time at each iteration. These objectives are captured by a hybrid reward function which is defined as:
\begin{equation}
	r \coloneqq  \alpha_1 r_{\text{feas}} + \alpha_2 r_{\text{gap}} - \alpha_3 r_{\text{time}},
\end{equation}
where \( \alpha_1, \alpha_2, \alpha_3 \geq 0 \) are weighting coefficients that balance the contributions of the feasibility reward \( r_{\text{feas}} \), the gap improvement reward \( r_{\text{gap}} \), and the time penalty \( -r_{\text{time}} \), respectively. The components of \( r \) are described below:

\begin{enumerate}
	\item \textbf{Feasibility reward:} Penalizes infeasible actions and rewards feasible ones:
	\begin{equation}
		r_{\text{feas}} \coloneqq
		\begin{cases}
			-\beta_1 & \text{if  the master problem is infeasible} \\ 
			\beta_2 & \text{if  the master problem is feasible}
		\end{cases},
	\end{equation}
	where \( \beta_1, \beta_2 \geq 0 \).
	
	\item \textbf{Gap improvement reward:} Encourages actions that reduce the GBD bound gap between iterations:
	\begin{equation}\label{eq:bound_improvement}
		r_{\text{gap}} \coloneqq
		\begin{cases}
			\displaystyle \left| \frac{(\text{UBD}_{{k-1}} - \text{LBD}_{{k-1}}) - (\text{UBD}_{{k}} - \text{LBD}_{{k}})}{\text{UBD}_{0} - \text{LBD}_{0}} \right| & \text{if  the master problem is feasible}, \\
			0 & \text{otherwise}
		\end{cases},
	\end{equation}
	where \( \text{UBD}_{{k}} \), \( \text{LBD}_{{k}} \) are the upper and lower bounds at the current iteration, \( \text{UBD}_{{k-1}} \), \( \text{LBD}_{{k-1}} \) denote their values from the previous iteration, and  \( \text{UBD}_{{0}} \), \( \text{LBD}_{{0}} \) are the initial bounds. This term only applies when the agent's action renders the master problem feasible; otherwise, the reward is zero since any bound improvement results from  the  MIP solver.
	
\item \textbf{Time penalty:}  Discourages actions that prolong the solution time of the subproblem. Given that the solution of the master problem can affect the complexity of the subproblem, this term allows the framework to incorporate the efficiency of the subproblem as one of its objectives. Let \( t_{\text{SP}} \) denote the time required to solve the subproblem. To prevent excessively large values from dominating $r$, the penalty is clipped at a threshold \( \tau \):  
\begin{equation}
	r_{\text{time}} \coloneqq \min(t_{\text{SP}}, \tau).
\end{equation}
\end{enumerate}

\subsubsection{Training Phase}
The agent is trained over a fixed number of episodes $N_{\text{ep}}$. Among the available actor-critic methods, we adopt the PPO algorithm~\cite{schulman2017proximal} due to its robustness and stability in complex training environments. In PPO, the critic estimates the state-value function, while the actor is updated based on an \emph{advantage} estimate, which measures how much better a selected action performs relative to the expected performance at the same state.
Agent training proceeds by collecting and storing experience tuples of the form $(S_k,\bm{a}_k,r,S_{k+1})$, where $S_k \coloneqq (A_{d_k}, X_{f_k}, X_{e_k})$ and $S_{k+1} \coloneqq (A_{d_{k+1}}, X_{f_{k+1}}, X_{e_{k+1}})$. These stored tuples are then used to update the actor $\pi_{\theta}$ and critic $V_{\phi}$ networks according to the PPO algorithm. The complete training loop for the agent within the GBD framework is summarized in Algorithm~\ref{alg:rl_gbd_step}.

\begin{algorithm}[!ht]
	\caption{Training procedure  of the RL agent.}
	\label{alg:rl_gbd_step}
	\small
	\DontPrintSemicolon
	\KwIn{\( \bm{y}^0 \), \( N_{\text{ep}} \), \( \varepsilon \), \( \pi_\theta \), \( V_\phi \), \( \text{LBD}_{\text{init}}, \text{UBD}_{\text{init}} \), \( \beta_1, \beta_2, \tau \), \( \alpha_1, \alpha_2, \alpha_3 \)}
	Initialize \( V_{\phi} \)\;
	\For{\( \text{episode} \gets 1 \) \KwTo \( N_{\text{ep}} \)}{
		Reset problem with new parameters\;
		Set \( \text{LBD}_{\text{cur}} \gets \text{LBD}_{\text{init}} \), \( \text{UBD}_{\text{cur}} \gets \text{UBD}_{\text{init}} \)\;
		Set \( \text{LBD}_{\text{prev}} \gets \text{LBD}_{\text{init}} \), \( \text{UBD}_{\text{prev}} \gets \text{UBD}_{\text{init}} \)\;
		Initialize \( \bm{y} \gets \bm{y}^0 \)\;
		Solve subproblem \( \mathcal{S}(\bm{y}) \)\;
		\eIf{subproblem \( \mathcal{S}(\bm{y}) \) is feasible}{
			Add $\mathcal{C}_o^{(k)}$  to master problem \( \mathcal{M} \)\;
			Update \( \text{UBD}_{\text{cur}} \gets \min(\text{UBD}_{\text{cur}}, Z(\bm{y})) \)\;
		}{
			Solve subproblem \( \mathcal{F}(\bm{y}) \)\;
			Add $\mathcal{C}_f^{(k)}$ to master problem \( \mathcal{M} \)\;
		}
		Encode  master problem \( \mathcal{M} \) as \( \mathcal{G}_t  = (A_{d_t}, X_{f_t}, X_{e_t})\)\;
		\While{\( \text{UBD}_{\text{cur}} - \text{LBD}_{\text{cur}} > \varepsilon \)}{
			
			Obtain  \( {\bm{a}}\) with \(\pi_{\theta} \text{ and } \mathcal{G}_t \) using Equation~\eqref{eq:action_selection} \;
			
			\eIf{\( {\bm{a}} \) is feasible for master problem \( \mathcal{M} \)}{
				Set \( r_{\text{feas}} \gets \beta_2 \)\;
				Evaluate  \(\mathcal{O}_k ({\bm{a}}) \) i.e., the RHS of $\mathcal{C}_o^{(k)}$ \( \forall  k \in K_O \)\;
				Update \( \text{LBD}_{\text{cur}} \gets  \max_{k \in K_O} \mathcal{O}_k ({\bm{a}})  \)\; 
				Compute \( r_{\text{gap}} \) using Equation~\eqref{eq:bound_improvement}\;
				Set \( \bm{y} \gets {\bm{a}} \)\;
			}{
				Set \( r_{\text{feas}} \gets -\beta_1,\ r_{\text{gap}} \gets 0 \)\;
				Solve  master problem \( \mathcal{M} \) to obtain \( \overline{\bm{y}}, \mu_b \)\;
				Set \( \bm{y} \gets \overline{\bm{y}} \), \( \text{LBD}_{\text{cur}} \gets \mu_b \)\;
			}
			Solve subproblem \( \mathcal{S}(\bm{y}) \)\;
			\eIf{subproblem \( \mathcal{S}(\bm{y}) \) is feasible}{
				Add $\mathcal{C}_o^{(k)}$ to master problem \( \mathcal{M} \)\;
				Update \( \text{UBD}_{\text{cur}} \gets \min(\text{UBD}_{\text{cur}}, Z(\bm{y})) \)\;
			}{
				Solve   subproblem \( \mathcal{F}(\bm{y}) \)\; Add $\mathcal{C}_f^{(k)}$ to  master problem \( \mathcal{M} \)\;
			}
			Measure  \( t_{\text{SP}} \); set \( r_{\text{time}} \gets \min(t_{\text{SP}}, \tau) \)\;
			Compute 
			\(
			r \gets \alpha_1 r_{\text{feas}} + \alpha_2 r_{\text{gap}} - \alpha_3 r_{\text{time}}
			\)\;
			Encode master problem \( \mathcal{M} \) as \( \mathcal{G}_{t+1}  = (A_{d_{t+1}}, X_{f_{t+1}}, X_{e_{t+1}})\)\;
			Store  \( (\mathcal{G}_t, {\bm{a}}, r, \mathcal{G}_{t+1}) \)\;
			Update  \( V_\phi \) and  \( \pi_\theta \) using PPO\;
			Set \( \mathcal{G}_t \gets \mathcal{G}_{t+1} \)\;
			Update \( \text{LBD}_{\text{prev}} \gets \text{LBD}_{\text{cur}} \), \( \text{UBD}_{\text{prev}} \gets \text{UBD}_{\text{cur}} \)\;
		}
	}
\end{algorithm}

\section{Verification and Coordination Mechanism  }\label{sec:proposed_approach_integration}
While the agent seeks to enhance the computational efficiency of GBD, relying solely on it to solve the master problem presents several challenges. First, although the agent’s reward structure encourages it to produce actions that are feasible with respect to the master problem, its actions are not guaranteed to be feasible. Therefore, a sole reliance on the agent could result in the passing of infeasible assignments to the subproblem. Second, the actions provided by the agent are not necessarily optimal. Consequently, poor actions by the agent could produce an LBD that exceeds the current best UBD. In such a case, the LBD resulting from the agent's action cannot be regarded as a valid lower bound for the original MINLP.

To address these challenges, we propose a verification and coordination strategy that leverages the agent to quickly generate assignments for the binary decision variables, with fallbacks to a MIP solver whenever feasibility is at risk or when the agent’s actions would cause the LBD to exceed the UBD. We refer to this approach as \textit{confidence-based assignment}. The main steps are summarized in Algorithm~\ref{alg:confidence_based_assignment}, and the details are provided in the next section.

\begin{algorithm}[!ht]
	\small
	\caption{Confidence-based assignment.}
	\label{alg:confidence_based_assignment}
	\DontPrintSemicolon
	\SetKwComment{Comment}{$\triangleright$ }{}
	
	\KwIn{\( \bm{p} \in [0,1]^m \), \( \delta_1, \delta_2 \) with \(0 \le \delta_1 \le \delta_2 \le 1\), \( \text{UBD} \),  \( \{\mathcal{O}_k(\cdot)\}_{k\in K_O} \)}
	\KwOut{\( \bm{y} \in \{0,1\}^m \), \( \mu_b \)}
	
	Initialize \( \bm{y} \) as free; \( \mathcal{Y} \gets \emptyset \)\;
	
	\ForEach{\( i \in \{1,\dots,m\} \)}{
		\eIf{\( p_i \le \delta_1 \)}{
			Set \( y_i \gets 0 \); add \( i \) to \( \mathcal{Y} \)\;
		}{
			\eIf{\( p_i \ge \delta_2 \)}{
				Set \( y_i \gets 1 \); add \( i \) to \( \mathcal{Y} \)\;
			}{
				Leave \( y_i \) free\;
			}
		}
	}
	
	\eIf{\( |\mathcal{Y}| = m \)}{ \Comment{Full assignment}
		\If{\( \bm{y} \) infeasible for master problem \( \mathcal{M} \)}{
			Solve master problem \( \mathcal{M} \) to obtain \( \overline{\bm{y}}, \overline{\mu}_b \)\;
			Set \( \bm{y} \gets \overline{\bm{y}},\ \mu_b \gets \overline{\mu}_b \)\;
		}
		\Else{
			\eIf{\( K_O \neq \emptyset \)}{
				Compute \( \hat{\mu}_b \gets \max_{k\in K_O} \mathcal{O}_k(\bm{y}) \)\;
			}{
				Set \( \hat{\mu}_b \gets -\infty \)\;
			}
			\eIf{\( \hat{\mu}_b \le \text{UBD} \)}{
				Set \( \mu_b \gets \hat{\mu}_b \)\;
			}{
				Solve master problem \( \mathcal{M} \) (with $y_1,\dots, y_m$ free) to obtain \( \overline{\bm{y}}, \overline{\mu}_b \)\;
				Set \( \bm{y} \gets \overline{\bm{y}},\ \mu_b \gets \overline{\mu}_b \)\;
			}
		}
	}{
		\eIf{\( 0 < |\mathcal{Y}| < m \)}{ \Comment{Partial assignment}
			Solve master problem \( \mathcal{M} \) \emph{with \( y_i \) fixed, \( \forall i \in \mathcal{Y} \)} to obtain \( \overline{\bm{y}}, \overline{\mu}_b \)\;
			\If{master problem \( \mathcal{M} \) infeasible}{
				Solve master problem \( \mathcal{M} \) (with $y_1,\dots, y_m$ free) to obtain \( \overline{\bm{y}}, \overline{\mu}_b \)\;
				Set \( \bm{y} \gets \overline{\bm{y}},\ \mu_b \gets \overline{\mu}_b \)\;
			}
			\Else{
				\eIf{\( \overline{\mu}_b \le \text{UBD} \)}{
					Set \( \bm{y} \gets \overline{\bm{y}},\ \mu_b \gets \overline{\mu}_b \)\;
				}{
					Solve master problem \( \mathcal{M} \) (with $y_1,\dots, y_m$ free) to obtain \( \overline{\bm{y}}, \overline{\mu}_b \)\;
					Set \( \bm{y} \gets \overline{\bm{y}},\ \mu_b \gets \overline{\mu}_b \)\;
				}
			}
		}{ \Comment{No assignment}
			Solve master problem \( \mathcal{M} \) (with $y_1,\dots, y_m$ free) to obtain \( \overline{\bm{y}}, \overline{\mu}_b \)\;
			Set \( \bm{y} \gets \overline{\bm{y}},\ \mu_b \gets \overline{\mu}_b \)\;
		}
	}
	
	\Return \( (\bm{y}, \mu_b) \)\;
\end{algorithm}

\subsection{Confidence-Based Assignments}\label{sec:confidence_based_assignment}
As discussed in Section~\ref{sec:action}, the agent's policy produces a probability vector $\bm{p} = (p_1,\dots,p_m)$, where each $p_i$ parameterizes a Bernoulli distribution for the corresponding binary decision variable $y_i$. During training, actions are sampled stochastically from these distributions to encourage exploration. However, during deployment, a deterministic strategy is adopted. To enhance reliability of the agent's actions, the probabilities are post-processed using a confidence-based rule.
Specifically, we define a high-confidence region as $[0,\delta_1] \cup [\delta_2,1]$ with $0 \leq \delta_1 \leq \delta_2 \leq 1$. The probabilities within this region are fixed deterministically to $0$ (if $p_i \leq \delta_1$) or $1$ (if $p_i \geq \delta_2$). The probabilities in the remaining interval $(\delta_1,\delta_2)$ are treated as uncertain, and the corresponding variables are left unfixed for the MIP solver to determine. This procedure yields three possible scenarios:
\begin{enumerate}
	\item \textbf{Full Assignment:} All probabilities lie within the confidence region, so all binary variables are fixed deterministically.
	\item \textbf{Partial Assignment:} A subset of variables is fixed deterministically, and the remaining ones are determined by the solver.
	\item \textbf{No Assignment:} None of the probabilities meet the confidence threshold, so all variables are determined by the solver.
\end{enumerate}

These are described in detail below.

\subsubsection{Full Assignment}
When all probabilities are within the high-confidence region, all binary variables are fixed to obtain a candidate assignment, denoted by \( \hat{\bm{y}} \). Two checks are performed to assess its suitability.
\begin{enumerate}
	\item \textbf{Feasibility:}  
	The candidate \( \hat{\bm{y}} \) is validated against the constraints of the master problem. This check is restricted to feasibility cuts and pure binary constraints. Specifically, for each feasibility cut \(\mathcal{C}_f^{(k)},~\forall k \in K_F\), and pure binary constraint \(\mathcal{C}_p\), the LHS is evaluated using \( \hat{\bm{y}} \). If any of these constraints are violated, \( \hat{\bm{y}} \) is deemed infeasible. In this case, the MIP solver is invoked to solve the master problem and the resulting solution is passed to the subproblem.
	
	\item \textbf{Evaluation of cost \( \hat{\mu}_b \):}  
	If \( \hat{\bm{y}} \) passes the feasibility check, the next step is to compute the cost \( \hat{\mu}_b \) associated with it. In this case, solving the master problem with \( \hat{\bm{y}} \) can be avoided. Instead, \( \hat{\mu}_b \) is obtained by evaluating the RHS of each optimality cut, $\mathcal{O}_k (\hat{\bm{y}}),~\forall k \in K_O$, and taking the maximum:
	\[
	\hat{\mu}_b \coloneqq \max_{k \in K_O} \mathcal{O}_k (\hat{\bm{y}}).
	\]
	Since the agent’s assignments are not guaranteed to be optimal, the successive values of \( \hat{\mu}_b \) may not exhibit the monotonic behavior of the LBD in classical BD. Thus, \( \hat{\mu}_b \) is treated as a \emph{candidate} LBD. The mechanism used to enforce monotonicity is discussed later in Section~\ref{sec:monotonicity}. To  assess the suitability of \( \hat{\bm{y}} \), \( \hat{\mu}_b \) is compared with the current best \( \mathrm{UBD} \). If \( \hat{\mu}_b > \mathrm{UBD} \), the candidate assignment is rejected, since it would yield an invalid lower bound. In this case, the solver is invoked to solve the master problem, and the resulting solution is passed to the subproblem. Otherwise, if \( \hat{\mu}_b \leq \mathrm{UBD} \), the candidate assignment is  forwarded to the subproblem. 
\end{enumerate}

\subsubsection{Partial Assignment}
When only a subset of probabilities lies within the high-confidence region, the corresponding variables are fixed, and the master problem is solved to determine the remaining ones. If the fixed values render the master problem infeasible, the assignment is discarded, and the solver is invoked to resolve the master problem with all binary variables free. The resulting solution is then passed to the subproblem.  

If the master problem is feasible for the fixed values, the solver returns a complete binary vector together with the corresponding cost \( \overline{\mu}_b \), which serves as a candidate LBD. To assess the suitability of this solution, \( \overline{\mu}_b \) is compared with the current best UBD. If \( \overline{\mu}_b \leq \mathrm{UBD} \), the solution is accepted and forwarded to the subproblem. Otherwise, the master problem is re-solved with all binary variables free, and the resulting solution is passed to the subproblem.

\subsubsection{No Assignment}
When none of the probabilities lies within the high-confidence region, the solver is invoked to solve the master problem directly. The resulting solution is then passed to the subproblem.

\begin{algorithm}[t]
	
	\caption{Graph-based agent integration into GBD.}
	\small
	\label{alg:learning_augmented_gbd}
	\DontPrintSemicolon
	\SetAlgoLined
	
	\KwIn{\( \bm{y}^0 \),  \( \varepsilon \),  \( \pi_{\theta} \),  \( \delta_1, \delta_2 \)}
	
	Initialize: \( k \gets 0 \), \( \text{LBD} \gets -\infty \), \( \text{UBD} \gets +\infty \), \( K_O \gets \emptyset \), \( K_F \gets \emptyset \)\;
	
	\While{\( \text{UBD} - \text{LBD} > \varepsilon \)}{
		Solve subproblem \( \mathcal{S}(\bm{y}^k) \)\;
		\eIf{subproblem \( \mathcal{S}(\bm{y}^k) \) is feasible}{
			Update \( \text{UBD} \gets \min(\text{UBD}, Z(\bm{y}^k)) \)\;
			Add $\mathcal{C}_o^{(k)}$ to  master problem \( \mathcal{M} \)\; 
			Update \( K_O \gets K_O \cup \{k\} \)\;
		}{
			Solve subproblem \( \mathcal{F}(\bm{y}^k) \)\;
			Add $\mathcal{C}_f^{(k)}$ to master problem \( \mathcal{M} \)\; 
			Update \( K_F \gets K_F \cup \{k\} \)\;
		}
		
		Construct \( \mathcal{G}_k = (A_{d_k}, X_{f_k}, X_{e_k}) \) of master problem \( \mathcal{M} \)\;
		Obtain  probabilities \( {\bm{p}} \gets \pi_{\theta}(A_{d_k}, X_{f_k}, L_k ) \)\;
		Invoke Algorithm~\ref{alg:confidence_based_assignment} with \( {\bm{p}}, \delta_1, \delta_2, \text{UBD}, \{\mathcal{O}_k (\bm{y})\}_{k \in K_O} \)   to obtain \( \bm{y}^{k+1}, \mu_b^k \)\;
		
		Update \( \text{LBD} \gets \max(\text{LBD}, \mu_b^k) \), \( \bm{y}^k \gets \bm{y}^{k+1} \), \( k \gets k + 1 \)\;
	}
\end{algorithm}

\subsection{Enforcing Monotonicity of the Lower Bound}\label{sec:monotonicity}
The cost $\mu_b$ resulting from Algorithm~\ref{alg:confidence_based_assignment} cannot be accepted directly as the new LBD, as it is not expected to exhibit the monotonic increase guaranteed in the classical GBD. To ensure monotonicity, the lower bound is explicitly updated at each iteration to enforce a nondecreasing sequence. Let $\mu_b$ denote the cost evaluated from the current master problem solution obtained from Algorithm~\ref{alg:confidence_based_assignment}, and let $\text{LBD}_{k-1}$ denote the lower bound from the previous iteration. The update rule used is
\[
\text{LBD}_k \coloneqq \max\left\{ \mu_b, \text{LBD}_{k-1} \right\}.
\]

\begin{remark}
	In the case of a full assignment, both the feasibility and cost-consistency checks in Algorithm~\ref{alg:confidence_based_assignment} are performed through direct evaluation of the feasibility cuts, pure binary constraints, and optimality cuts, without solving any additional optimization problem. 
	For partial assignments, feasibility can sometimes be verified directly. However, depending on the structure of the free variables and their interaction with the constraints, checking feasibility may require solving a reduced master problem. 
	Similarly, cost consistency in the partial assignment case requires solving a reduced master problem. 
	Although this introduces some additional effort, a reduced master problem is  much smaller than the full master problem. 
\end{remark}

\section{Case Studies}\label{sec:simulations}
Two case studies are presented to demonstrate the performance of the proposed approach. The first is a generic MINLP that provides a controlled setting to evaluate the specific design choices of the proposed method. The second focuses on a real application, where the proposed approach is applied to an irrigation scheduling problem.

\subsection{Case Study 1}
First, we consider the following problem adapted from~\cite{floudas1995nonlinear} denoted as $\mathcal{E}$:
\begin{align}
	\min_{\bm{x},\bm{y}} \quad &
	5y_{1} + 8y_{2} + 6y_{3} + 10y_{4} + 6y_{5} 
	- 10x_{3} - 15x_{5} - 15x_{9} + 15x_{11} + 5x_{13} - 20x_{16} \nonumber\\
	& + \exp(x_{3}) + \exp\!\left(\frac{x_{5}}{1.2}\right) 
	- 60\ln(x_{11} + x_{13} + 1) + 140 \tag{E1} \label{E1} \\
	\text{s.t.} \quad
	& -\ln(x_{11} + x_{13} + 1) \le 0 \tag{E2} \label{E2} \\
	& - x_{3} - x_{5} - 2x_{9} + x_{11} + 2x_{16} \le 0 \tag{E3} \label{E3} \\
	& - x_{3} - x_{5} - 0.75x_{9} + x_{11} + 2x_{16} \le 0 \tag{E4} \label{E4} \\
	& x_{9} - x_{16} \le 0 \tag{E5} \label{E5} \\
	& 2x_{9} - x_{11} - 2x_{16} \le 0 \tag{E6} \label{E6} \\
	& -0.5x_{11} + x_{13} \le 0 \tag{E7} \label{E7} \\
	& 0.2x_{11} - x_{13} \le 0 \tag{E8} \label{E8} \\
	& \exp(x_{3}) - U y_{1} \le 1 \tag{E9} \label{E9} \\
	& \exp\!\left(\frac{x_{5}}{1.2}\right) - U y_{2} \le 1 \tag{E10} \label{E10} \\
	& 1.25x_{9} - U y_{3} \le 0 \tag{E11} \label{E11} \\
	& x_{11} + x_{13} - U y_{4} \le 0 \tag{E12} \label{E12} \\
	& -2x_{9} + 2x_{16} - U y_{5} \le 0 \tag{E13} \label{E13} \\
	& y_{1} + y_{2} = 1 \tag{E14} \label{E14} \\
	& y_{4} + y_{5} \le 1 \tag{E15} \label{E15} \\
	& \bm{y}^{T} = (y_{1},y_{2},y_{3},y_{4},y_{5})  \in \{0,1\}^{5}, \quad 
	\bm{x}^{T}  = (x_{3},x_{5},x_{9},x_{11},x_{13},x_{16}) \in \mathbb{R}^{6} \nonumber\\
	& \bm{a}^{T} = (0,0,0,0,0,0), \quad 
	\bm{b}^{T} = (2,2,2,\text{--},\text{--},3), \quad 
	\bm{a} \le \bm{x} \le \bm{b}, \quad U = 10 \nonumber
\end{align}
\normalsize

\subsubsection{Simulation Settings}
\paragraph{Data Generation}\label{sec:data_generation_ex_1}
Several components of problem $\mathcal{E}$ were parameterized to produce varied solutions across simulation runs. In particular, the cost function was expressed as
\[
c_1y_{1} + c_2y_{2} + c_3y_{3} + c_4y_{4} + c_5y_{5} + F_{\mathcal{E}}(\bm{x}),
\]
where $F_{\mathcal{E}}(\bm{x})$ compactly denotes all $\bm{x}$-dependent terms in Equation~\eqref{E1}. The coefficients $c_1, c_2, c_3, c_4$ were sampled as positive integers from the interval $[1,39]$, while $c_5$ was also sampled as a positive integer from $[1,7]$. These ranges were selected based on a series of simulation experiments, which indicated that they yield parameter realizations with varied optimal solutions.

Before training the agent, 3{,}000 distinct realizations of the coefficient vector $(c_1,\dots,c_5)$ were generated. For each parameter set, GBD was used to solve the corresponding instance. The specific master problem and subproblem formulations used during this process are provided in  Section 1.1 of the Supplementary Information. 
The master problem was solved using Gurobi~12.0.1~\cite{gurobi}, while the subproblem was solved using IPOPT~3.12.6~\cite{wachter2006ipopt}.

\paragraph{Agent Training}
The dataset from Section~\ref{sec:data_generation_ex_1} was used to train the agent in the imitation learning stage. The model architecture and hyperparameters are summarized in Table 2 of the Supplementary Material. Graph operations were implemented using the Spektral library~\cite{grattarola2020spektral}, and training was carried out in Keras~\cite{chollet2015keras}. In the RL stage, training was performed over 10,000 episodes, each corresponding to a distinct parameter realization. The critic network architecture, reward design parameters, and RL training hyperparameters are presented in Table~3 of the Supplementary Material. The training of the model in the imitation learning and RL stages took 20 minutes and 3 hours, respectively on an Intel\textregistered\ Core\texttrademark\ i7-14700 processor.

\paragraph{Performance Evaluation}
In addition to the proposed approach, two other approaches were considered: (i) the IL agent and (ii) an RL agent with a randomly initialized actor network. For evaluation, a comparative study was conducted using 100 new parameter realizations. In this study, each of the three agent-based methods, together with the proposed approach, was embedded into GBD following the Algorithm~\ref{alg:learning_augmented_gbd}. For reference, classical GBD was also applied to the same 100 problem instances. The confidence thresholds were set at $\delta_1 = 0.10$ and $\delta_2 = 0.90$. These values were manually tuned to strike a balance between computational efficiency and the selection of assignments with low uncertainty. Performance was assessed on the basis of runtime, evolution of the Benders gap, and convergence to the true solution. In addition, for agent-based approaches, the frequency with which their assignments were accepted was also recorded.

\subsubsection{Results and Discussion}
\begin{figure*}[t]
	\centerline{\includegraphics[width=\textwidth]{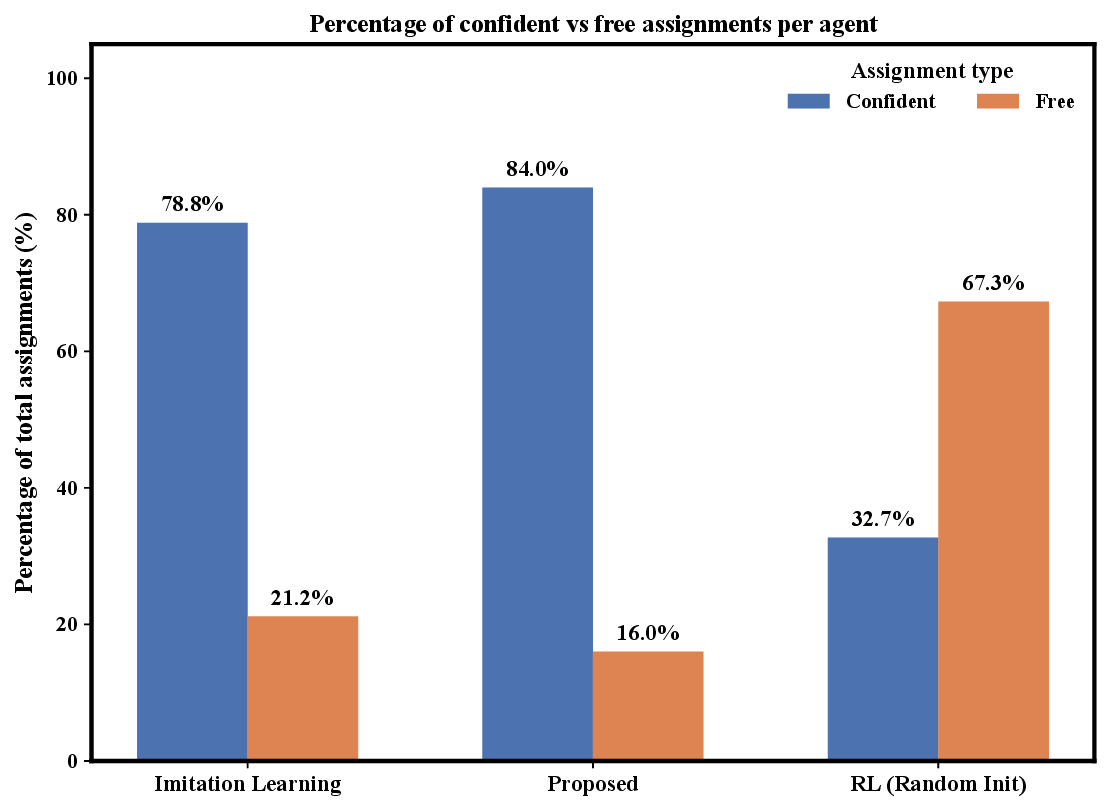}}
	\caption{Percentage of confident and free binary variable assignments produced by each agent-based method over 100 test instances.} 
	\label{fig:assignment_frequency}
\end{figure*}

Table~\ref{tbl:runtime_metrics} summarizes the runtime performance of the four approaches. For this example, the master problem dominates the total runtime of the GBD, accounting for approximately 84\% of the total time. Both the proposed approach and the IL agent reduced the runtime of the master problem relative to classical GBD, with reductions of 56\% and 42\%, respectively. By contrast, the RL agent with a randomly initialized actor increased master problem runtime by 4\%, indicating that its assignments often prolonged the solution process. In terms of total GBD runtime, the proposed method achieved a reduction 42\%, the IL agent achieved a 31\% reduction, while the randomly initialized RL agent increased total time by 3\%. These results are consistent with the trends in Figure~\ref{fig:assignment_frequency}, which shows that the proposed method determined binary variables 84\% in 100 instances, compared to 78.8\% for the IL agent and 32. 7\% for the RL agent initialized randomly.

These findings confirm our hypothesis that the combination of IL and RL yields greater computational benefits in combinatorial settings than the use of IL alone or RL alone. In the proposed approach, the RL stage refines the IL policy, producing more confident and higher-quality binary assignments, as reflected in Figure~\ref{fig:assignment_frequency}. Consequently, while both the proposed method and the IL agent reduce the master problem runtime (by 56\% and 42\%, respectively) and the total GBD runtime (by 42\% and 31\%, respectively), the proposed approach achieves greater reductions, highlighting the benefits of the two-stage design. By contrast, the randomly initialized RL agent struggles in the complex GBD environment, producing fewer suitable assignments and triggering more frequent solver calls. This limitation is also reflected in the convergence results of Table~\ref{tbl:convergence_metrics}, where the randomly initialized RL agent fails to reach the optimal solution in 30 of 100 instances, with relative optimality gaps between 0.09 and 0.25. In comparison, both the IL agent and the proposed method converge in all instances.

\begin{figure*}[t]
	\centerline{\includegraphics[width=\textwidth]{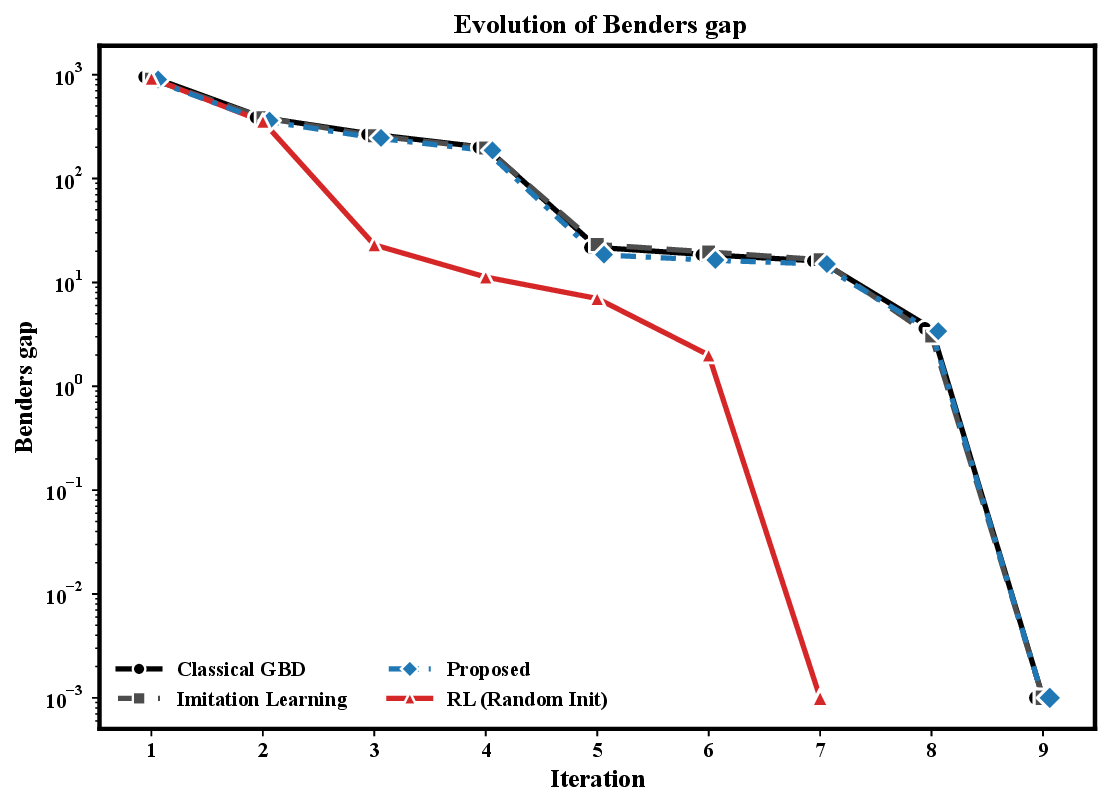}}
	\caption{Evolution of the median Benders gap across 20 problem instances, plotted over the median number of iterations for each solution method.} 
	\label{fig:benders_gap}
\end{figure*}

Figure~\ref{fig:benders_gap} shows the evolution of the Benders gap, evaluated over the median number of iterations for each method (7 for the randomly initialized RL agent and 9 for the other three). For each method, 20 instances terminating at the median iteration count were selected and the median gap between these instances was calculated at each iteration. Although the randomly initialized RL agent terminated in fewer iterations, it produced more computationally expensive master problems (Table~\ref{tbl:runtime_metrics}), eliminating any potential runtime benefits. Importantly, both the imitation learning agent and the proposed method exhibited convergence profiles closely aligned with the classical GBD.  
\begin{table}[t]
	\centering
	\footnotesize
\caption{Runtime performance across methods over 100 test instances. Values are reported as mean $\pm$ standard deviation over the 100 instances. Percentages in parentheses indicate improvement relative to classical GBD.}
	\label{tbl:runtime_metrics}
	\begin{tabular}{lcc}
		\toprule
		\textbf{Method} & \textbf{Master problem runtime (s)} & \textbf{Total GBD runtime (s)}  \\
		\midrule
		Classical GBD     & $0.54 \pm 0.13$ & $0.64 \pm 0.13$ \\
		Imitation learning Agent          & $0.31 \pm 0.12$ (42\%) & $0.44 \pm 0.12$ (31\%) \\
		\textbf{Proposed} & $\mathbf{0.24 \pm 0.11}$ (\textbf{56\%}) & $\mathbf{0.37 \pm 0.12}$ (\textbf{42\%}) \\
		RL (Random Init)  & $0.56 \pm 0.31$ (-4\%) & $0.66 \pm 0.33$ (-3\%)  \\
		\bottomrule
	\end{tabular}
\end{table}
\begin{table}[t]
	\centering
	\footnotesize
\caption{Convergence performance across methods over 100 test instances. Relative cost gap values are reported as mean $\pm$ standard deviation over the 100 instances.}
	\label{tbl:convergence_metrics}
	\begin{tabular}{lcc}
		\toprule
		\textbf{Method} & \textbf{Converged/Total} & \textbf{Relative gap} \\
		\midrule
		Classical GBD     & 100/100 & $0.00 \pm 0.00$ \\
		Imitation learning Agent          & 100/100 & $0.00 \pm 0.00$ \\
		Proposed          & 100/100 & $0.00 \pm 0.00$ \\
		RL (Random Init)  &  70/100 & $0.09 \pm 0.16$ \\
		\bottomrule
	\end{tabular}
\end{table}
\normalsize

\subsection{Case Study 2: Closed-loop Irrigation Scheduling}
In this case study, the proposed approach is applied to a closed-loop irrigation scheduling problem formulated as a mixed-integer MPC (MIMPC). The MIMPC scheduler was developed in~\cite{agyeman2024learning}, where it demonstrated improvements in water use efficiency compared to conventional scheduling methods. Despite these benefits, computational efficiency remained a key bottleneck. Since the scheduler must be solved on a daily basis to account for weather forecast uncertainty and because scaling to larger fields further increases complexity, decomposition methods such as GBD were recommended as a more efficient solution method. We apply the proposed approach to the MIMPC irrigation scheduler and evaluate its potential to improve computational efficiency without compromising its performance.

In irrigation scheduling, the goal is to determine when to irrigate within a scheduling horizon and how much water to apply on those days to ensure optimal crop development. Formally, the scheduling problem is stated as follows:

\noindent\textbf{Given:}
\begin{itemize}
	\item \textbf{Scheduling Horizon ($N$)}: Number of days over which irrigation decisions are to be made.
	\item \textbf{Weather Predictions}: Daily forecasts of reference evapotranspiration ($\text{ET}_0$) and precipitation ($\text{R}_{\text{n}}$), where $\text{ET}_0$ quantifies combined water loss due to evaporation and plant transpiration, and $\text{R}_{\text{n}}$ represents natural water input into the soil.
	\item \textbf{Crop Information}: Crop coefficient ($\text{K}_{\text{c}}$) values used to scale $\text{ET}_0$ to reflect crop-specific water requirements at different stages of growth.
	\item \textbf{Soil Moisture Content}: Initial spatial distribution of soil moisture ($\theta$) in the field.
\end{itemize}

\noindent\textbf{Determine:}
\begin{itemize}
	\item \textbf{Irrigation Timing}: The days within the scheduling horizon on which the irrigation is to be performed.
	\item \textbf{Irrigation Amounts}: Daily water application quantities ($u$) during each irrigation event, which is constrained to be zero on non-irrigation days and strictly positive on irrigation days.
\end{itemize}

In the MIMPC formulation, the irrigation timing is modeled with binary decision variables. Specifically, for each day on the horizon, a binary variable $c_j$ is introduced which is equal to one if irrigation is performed on day $j$ and zero otherwise.
The MIMPC scheduler employs a predictive soil moisture model to generate irrigation schedules that maintain soil moisture $\theta$ within a zone that is known to promote optimal crop development. In addition, the formulation minimizes the operating costs of irrigation events and the total volume of water applied. For day $d$, the scheduler $\mathcal{I}(\theta)$ is given by:
\begin{subequations}
	\small
	\begin{align}
		\min_{\bm{\bar{\epsilon}},~\bm{\underline{\epsilon}},~\bm{u},~\bm{c}} 
		& \quad 
		\underbrace{\sum_{j=d+1}^{d+N} \left[ \bar{Q} \bar{\epsilon}^2_j + \underline{Q} \underline{\epsilon}^2_j \right]}_{\text{Zone tracking}}
		+ \underbrace{\sum_{j=d}^{d+N-1} R_c c_j}_{\text{Irrigation cost}} 
		+ \underbrace{\sum_{j=d}^{d+N-1} R_u u_j}_{\text{Irrigation volume cost}}
		\label{eq:obj} \\[0.5em]
		\text{s.t.} \quad 
		& \theta_{j+1} = \mathcal{F}_m\left(\theta_j, K_{\text{c}_j}, \text{ET}_{0_j}, \text{R}_{\text{n}_j} u_j\right),
		&& j \in [d,\, d+N-1] \label{eq:cons1} \\[0.25em]
		& \theta_d = \theta(d), && \label{eq:cons2} \\[0.25em]
		& \underline{\nu} - \underline{\epsilon}_j \leq \theta_j \leq \bar{\nu} + \bar{\epsilon}_j, 
		&& j \in [d+1,\, d+N] \label{eq:cons3} \\[0.25em]
		& c_j\,\underline{u} \leq u_j \leq c_j\,\bar{u}, 
		&& j \in [d,\, d+N-1] \label{eq:cons4} \\[0.25em]
		& c_j \in \{0, 1\}, && j \in [d,\, d+N-1] \label{eq:cons5} \\[0.25em]
		& \underline{\epsilon}_j \geq 0, \quad \bar{\epsilon}_j \geq 0, 
		&& j \in [d+1,\, d+N] \label{eq:cons6}
	\end{align}
\end{subequations}
\normalsize
where $\bm{\bar{\epsilon}}\coloneqq [ \bar{\epsilon}_{d+1}, \dots, \bar{\epsilon}_{d+N}]$, $\bm{\underline{\epsilon}}\coloneqq [ \underline{\epsilon}_{d+1}, \dots, \underline{\epsilon}_{d+N}]$, $\bm{c}\coloneqq [ c_{d}, \dots, c_{d+N-1}]$, and $\bm{u}\coloneqq [u_{d}, \dots, u_{d+N-1}]$. The variables $\underline{\epsilon}_j$ and $\bar{\epsilon}_j$ are nonnegative slack variables introduced to relax, respectively, the lower $\underline{\nu}$ and upper $\bar{\nu}$ bounds of the target zone, with $\underline{Q}$ and $\bar{Q}$ denoting their per-unit penalty costs. The parameters $R_u$ and $R_c$ correspond to the per-unit cost of the applied irrigation volume and the fixed operational cost of the irrigation system, respectively.

The constraint~\eqref{eq:cons2} specifies the initial state, which is assumed to be known and obtained from the convergent soil moisture estimates produced in an offline estimation step.  
The constraint~\eqref{eq:cons1} represents the dynamics of soil moisture, initially described by a well-calibrated mechanistic model based on the Richards equation. The calibration uses field-specific parameters determined through an offline parameter estimation procedure. Due to the high computational burden of the mechanistic model, a long short-term memory (LSTM) network, which was trained according to the approach proposed in~\cite{agyeman2024learning}, is adopted as a surrogate to approximate the calibrated dynamics.  
The constraint~\eqref{eq:cons4} enforces operational limits on irrigation amounts such that if $c_j = 0$, the applied volume is fixed at zero. Conversely, if $c_j = 1$, the applied amount is constrained between $\underline{u}$ and $\bar{u}$.

The MIMPC scheduler was applied to a section of a large-scale agricultural field at the Lethbridge College Research Farm, located in southern Alberta, Canada. The field is equipped with a five-span center pivot irrigation system to enable precise and efficient water application. Before the application of the scheduler, field data collected from May 16, 2022, to August 31, 2022, were used to estimate both soil moisture and soil hydraulic parameters using the extended Kalman filtering approach. Moisture estimates were used to initialize the scheduler, while hydraulic parameters were used to calibrate the Richards equation describing soil moisture dynamics.

\subsubsection{Simulation Settings}
\paragraph{Data Generation}~\label{sec:data_generation_ex_2}
The data generation process involved solving the scheduler in a range of weather conditions. The weather inputs included forecasts of reference evapotranspiration, precipitation, and daily average temperature. Temperature forecasts, together with the empirical equations reported in~\cite{agyeman2024learning}, were used to compute crop coefficients for wheat, the crop cultivated in the study area. Historical weather data for the growing seasons 2009 to 2023 were obtained from the Alberta Weather Station database~\cite{alberta_weather_data}. For each year, data covering the wheat growing period (5~May to 4~September) were extracted. For each growing season, classical GBD was applied to determine irrigation schedules over the specified period, using a scheduling horizon of 14~days. This implementation followed a receding horizon approach, where the calibrated Richards equation was used to simulate field conditions. The specific formulations of the master problem and the subproblems solved during GBD are provided in Section 1.2 of the Supplementary Material, and the relevant parameters of $\mathcal{I}(\theta)$ used during data generation are listed in Table 1 of the Supplementary Material. The master problem was solved using Bonmin~1.8.9~\cite{bonami2008bonmin}, while the subproblem was solved using IPOPT~3.12.6~\cite{wachter2006ipopt}.

\paragraph{Agent Training}
The dataset from Section~\ref{sec:data_generation_ex_2} was used to train the agent in the IL stage. For this example, the network architecture and hyperparameters listed in Table 2 of the Supplementary Material were adopted, the only change being the output layer. Given the 14-day scheduling horizon of the scheduler, the output layer consisted of 14 sigmoid units. In the RL stage, the adopted critic architecture, the reward design parameters, and the hyperparameters of RL training are listed in Table 3 of the Supplementary Information. To promote generalization, RL training incorporated historical weather data from two nearby locations (Edmonton and Three Hills Weather Stations) and was carried out over 15,000 episodes. The agent was trained on an Intel\textregistered\ Core\texttrademark\ i7-14700 processor, where the training in the IL stage took 20 minutes while the training in the RL stage took 72 hours.

\paragraph{Performance Evaluation}
The proposed approach was used to generate wheat irrigation schedules during the 2024 growing season. The confidence thresholds were set at $\delta_1 = 0.10$ and $\delta_2 = 0.90$ to balance computational efficiency with assignment certainty. For comparison, the classical GBD was also applied.
A receding-horizon implementation was adopted, in which the calibrated Richards equation was used to represent the actual field dynamics. The simulation used observed 2024 weather data for actual field conditions, while noisy versions of weather inputs were used in scheduler evaluations.

Performance was compared in: (i) predicted crop yield, calculated using the method in~\cite{agyeman2024learning}; (ii) total prescribed irrigation; and (iii) efficiency of irrigation water use (IWUE), which is defined as predicted yield divided by total prescribed irrigation. Computational performance was analyzed over the 123 daily evaluations conducted from May 5 to September 4, 2024. The assessment included the daily master problem evaluation time,  the number of GBD iterations, the daily subproblem evaluation time, the daily scheduler evaluation time, and the total simulation time for the entire growing season.

\subsubsection{Results and Discussion}
\paragraph{Computational Performance}
Table~\ref{tbl:time_metrics_irrig} summarizes the computational performance of the proposed approach relative to classical GBD.  The relatively high standard deviations in the runtimes reflect the day-to-day variability of the scheduling problem, which arises from differences in weather conditions and soil moisture states throughout the growing season.
The proposed method reduced the runtime of the master problem by 38\%, the runtime on a daily basis by 23\%, and the total runtime on a seasonal basis by 23\%. However, for the irrigation scheduling problem, the subproblem was the primary bottleneck, since it embeds soil moisture dynamics modeled with an LSTM surrogate of the Richards equation. Although simplified, the surrogate retains sufficient complexity, contributing to the high subproblem runtimes. Importantly, accelerating the master problem also improved the efficiency of the subproblem, resulting in a reduction 21\% in the runtime of the subproblem. Two factors explain this effect. First, the proposed approach converged in slightly fewer iterations on average (10 versus 11 for classical GBD), reducing the total number of subproblem evaluations. Second, even after normalizing by iteration count, subproblem evaluations were faster under the proposed method (97~s per subproblem versus 112~s for classical GBD). This provides direct evidence of the agent’s contribution, as its reward structure explicitly penalizes long subproblem times and encourages assignments that yield more tractable subproblems. Together, these findings show that accelerating the master problem can translate into meaningful efficiency gains for the subproblem as well.

\begin{figure*}[t]
	\centerline{\includegraphics[width=\textwidth]{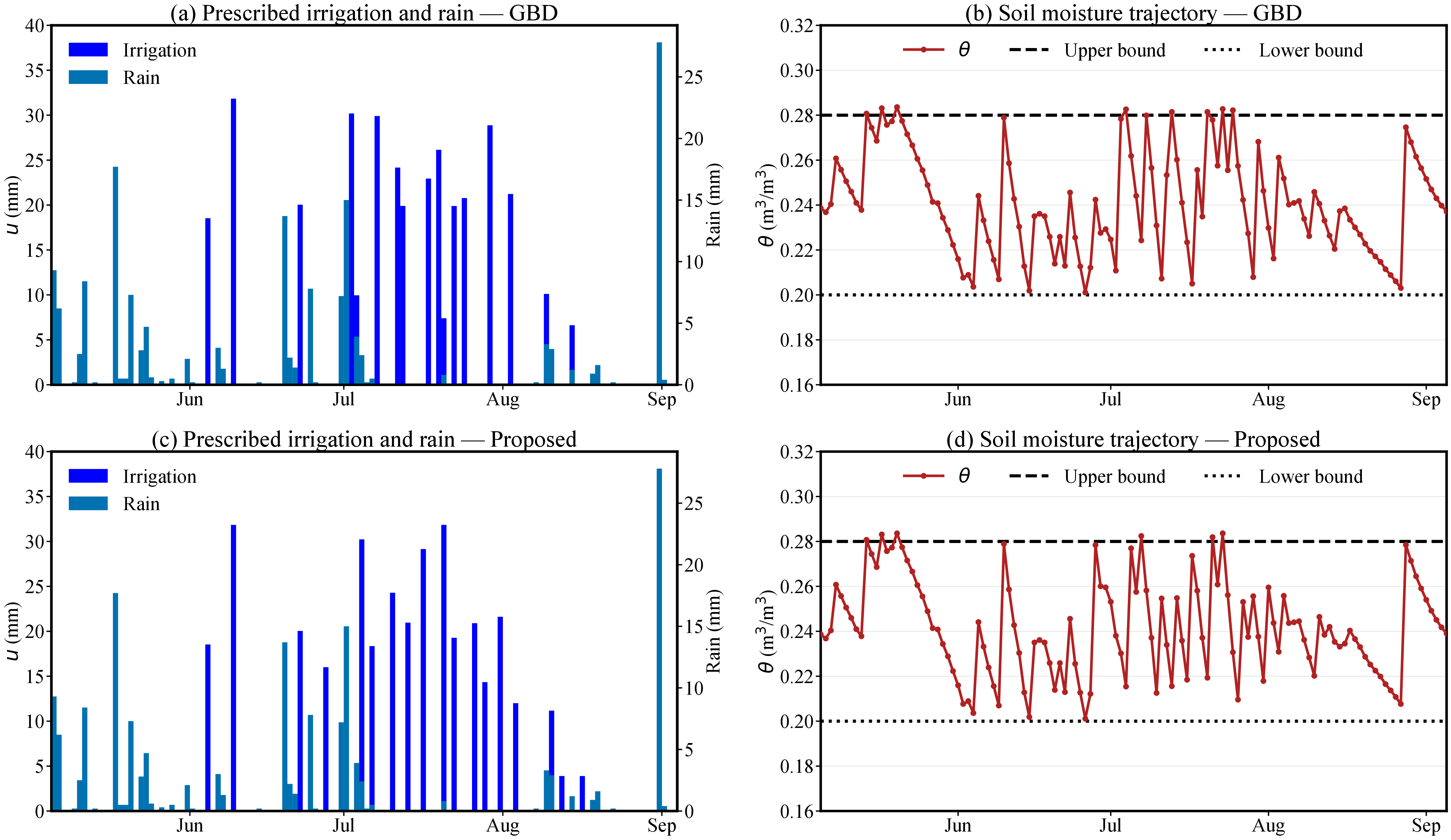}}	
		\caption{Comparison of prescribed irrigation schedules and corresponding soil moisture trajectories for the 2024 growing season.} 
	\label{fig:irrig_results}
\end{figure*}

\paragraph{Irrigation Performance}
Figure~\ref{fig:irrig_results} presents the prescribed irrigation schedules and soil moisture trajectories for the 2024 growing season under both the classical GBD and the proposed approach. In both cases, soil moisture was maintained successfully within the target zone, and the total prescribed seasonal irrigation was identical (Table~\ref{tbl:perf_metrics}). However, Figures~\ref{fig:irrig_results}(a) and (c) show differences in the timing of irrigation events. Except for the first three events, which occur on the same days in both methods, subsequent events are scheduled differently. This discrepancy is expected and is due to the hybrid nature adopted during the agent deployment. In some MPC evaluations, the agent fixes a subset of binary irrigation decisions and the solver then optimizes the remaining variables. This changes the decision space available to the solver. Since MPC minimizes cost over the prediction horizon, the solver adapts the unfixed decisions to best meet the objective. Consequently, the daily irrigation schedules from the proposed method can differ from those of the classical GBD, where the solver retains full control over all binary irrigation decisions.
Furthermore, Figures~\ref{fig:irrig_results}(b) and (d) show that the proposed approach led to fewer and smaller violations of the upper bound of the target zone. This reduction is associated with a 5\% increase in predicted crop yield. Although both approaches applied the same total irrigation volume, the proposed method used this water more effectively, resulting in a 5\% improvement in IWUE.

\begin{table}[t]
	\centering
	\footnotesize
	\caption{Comparison of computational runtimes between the proposed and classical GBD approaches over the entire growing season (5~May to 4~September~2024). Values are reported as mean~$\pm$~standard deviation across the 123 daily evaluations.}
	\begin{tabular}{lccc}
		\toprule
		\textbf{Metric} & \textbf{Proposed} & \textbf{Classical GBD} & \textbf{Improvement (\%)} \\
		\midrule
		Master problem runtime (s) & \textbf{83~$\pm$~18} & 133~$\pm$~26 & \textbf{38} \\
		Subproblem runtime (s) & \textbf{973~$\pm$~1192} & 1233~$\pm$~1461 & \textbf{21} \\
		Daily runtime (s) & \textbf{1056~$\pm$~1208} & 1366~$\pm$~1484 & \textbf{23} \\
		GBD iterations (count) & \textbf{10~$\pm$~12} & 11~$\pm$~12 & \textbf{9} \\
		\midrule
		Total evaluation time (min) & \textbf{2165} & 2800 & \textbf{23} \\
		\bottomrule
	\end{tabular}
	\label{tbl:time_metrics_irrig}
\end{table}

\begin{table}[t]
	\centering
	\footnotesize
	\caption{Irrigation performance comparison for the 2024 growing season.}
	\begin{tabular}{lccc}
		\toprule
		\textbf{Metric} & \textbf{Proposed} & \textbf{Classical GBD} & \textbf{Improvement (\%)} \\
		\midrule
		Total irrigation (m) & 0.35 & 0.35 & -- \\
		Predicted yield (kg/m$^2$) & \textbf{0.84} & 0.80 & \textbf{5} \\
		IWUE (kg/m$^3$) & \textbf{2.41} & 2.30 & \textbf{5} \\
		\bottomrule
	\end{tabular}
	\label{tbl:perf_metrics}
\end{table}

\subsection{Discussion}

The results presented in this work show that the graph-based   agent achieves a meaningful acceleration of the master problem. In addition, the results of the irrigation scheduling example demonstrate that the graph-based agent can also enhance the computational performance of the subproblem. This arises because, while the proposed approach focuses primarily on the master problem, it also considers as a secondary objective the computational tractability of the subproblem that results from the agent’s actions. Consequently, the approach proves beneficial in applications where the master problem constitutes the primary computational bottleneck (as in Case Study 1), as well as in cases where the subproblem dominates the computational effort (as in the irrigation scheduling case study). In general, these findings highlight that the proposed method improves the computational efficiency of GBD as a whole.

Several aspects of the proposed approach require some discussion. One limitation of the confidence-based assignment strategy in Algorithm~\ref{alg:confidence_based_assignment} concerns the certification of candidate LBDs. The filtering step, which compares each candidate LBD with the current best UBD, removes assignments that would yield invalid LBDs. However, it does not guarantee that all accepted assignments will yield valid lower bounds of the original MINLP. Although this limitation did not affect the quality of the results reported here, it highlights the need for stronger certification mechanisms in future studies. A stronger certification could be achieved by running the master problem under a limited computational budget, using the agent’s solution as a warm start to provide a valid lower bound. Although such measures would introduce some overhead, they could strengthen the theoretical guarantees of the proposed approach. 

A second issue is the computational effort for the RL stage. In the first case study, training took about 3~hours, while in the irrigation scheduling case, it required approximately 72~hours on an Intel\textregistered\ Core\texttrademark\ i7-14700 processor. Although these requirements represent a substantial initial investment, training is performed offline. Once trained, the agent can be deployed at a negligible cost.

Another issue concerns the choice of solvers. Although the IL stage requires data generated from an existing solver (Gurobi in Example~1 and Bonmin in Example~2), the framework itself is not tied to a particular solver. In practice, other solvers could be used to generate data, and the RL stage adapts the policy beyond the solver-specific behavior. Consequently, the computational gains observed here, particularly in percentage terms, are expected to be reproducible with other modern solvers.

Future work may also explore incorporating constraints directly into the agent design. Advances in constrained RL, such as constrained policy optimization, could be helpful in this context. Finally, while this work has focused on accelerating the master problem, a natural extension is to consider the subproblem. Unlike binary variables in MP, subproblem involves continuous decision variables and quantifying their prediction confidence  is more challenging. However, recent advances, such as Bayesian neural networks~\cite{kwon2020uncertainty} and ensemble-based approaches~\cite{rahaman2021uncertainty}, offer promising tools to assess the reliability of continuous-valued predictions.

\section{Conclusions}\label{sec:conclusion}
This work has demonstrated that graph-based learning can effectively accelerate the master problem in GBD. The proposed framework introduces an agent that processes a bipartite graph of the master problem to predict the values of its binary decision variables. The agent is trained using a two-stage design that integrates IL and RL. In the first stage, a graph neural network with a multi-headed output layer is trained via behavioral cloning. In the second stage, the policy from the first stage is refined within an actor–critic framework using a reward that balances feasibility and computational efficiency. Integration into GBD is achieved through a confidence-based assignment strategy, which leverages the agent for high-confidence predictions while falling back to a solver whenever feasibility is at risk or when assignments yield invalid lower bounds. 

The effectiveness of the approach was illustrated through two case studies. In a generic MINLP, the method accelerated the master problem without compromising the quality of the solution. The results further showed that IL provides a meaningful initial acceleration, while RL builds on this to deliver additional savings without compromising convergence to optimal solutions. In a closed-loop irrigation scheduling problem, formulated as a mixed-integer MPC, the framework accelerated GBD while preserving irrigation performance.

\section*{Acknowledgements}

The partial financial support of NSF CBET (award number 2313289) is gratefully acknowledged. IM would like to acknowledge financial support from the McKetta Department of Chemical Engineering.

\section*{Data Availability}
The numerical data for Figures~\ref{fig:assignment_frequency} and~\ref{fig:benders_gap} are provided in the Supplementary Material archive under \texttt{data/processed/} as \texttt{fig4\_assignment\_frequencies.xlsx} and \texttt{fig5\_gap\_curves.xlsx}, respectively. The time-series data for Figure~\ref{fig:irrig_results} are included in \texttt{data/processed} as \texttt{fig6\_data.xlsx}. Plotting scripts  are available in \texttt{code/}  as \texttt{fig4\_plot.py}, \texttt{fig5\_plot.py},  and \texttt{fig6\_plot.py}.

\section*{Code Availability}
The code for all case studies is available at \url{https://github.com/btagyeman/graph_il_rl_for_benders}.

\bibliographystyle{ieeetr}
\bibliography{references}
\end{document}